\documentclass[%
aip,
amsmath,amssymb,
preprint,%
]{revtex4-1}
\usepackage[dvipsnames]{xcolor}
\usepackage{graphicx}
\usepackage{dcolumn}
\usepackage{bm}
\usepackage{xr}
\usepackage{mathrsfs}
\usepackage[pdfborder={0 0 0},colorlinks=true,citecolor=blue,linkcolor=blue,CJKbookmarks=true,citecolor=blue]{hyperref}
\usepackage[utf8]{inputenc}
\usepackage[T1]{fontenc}
\usepackage{mathptmx}
\usepackage{etoolbox}
\usepackage{subfigure}
\usepackage[ruled,lined]{algorithm2e}

\newtheorem{theorem}{Theorem}
\newtheorem{lemma}[theorem]{Lemma}

\newtheorem{remark}[theorem]{Remark}
\makeatletter
\def\@email#1#2{%
	\endgroup
	\patchcmd{\titleblock@produce}
	{\frontmatter@RRAPformat}
	{\frontmatter@RRAPformat{\produce@RRAP{*#1\href{mailto:#2}{#2}}}\frontmatter@RRAPformat}
	{}{}
}%
\makeatother

\begin{document}

\title{ A new approach for stability analysis of $1-$D wave equation with time delay}
\author{Shijie Zhou}	
\affiliation{Research Institute of Intelligent Complex Systems, Fudan University, Shanghai 200433, China}
\author{Hongyinping Feng}
\affiliation{School of Mathematical Sciences, Shanxi University}	
\author{Zhiqiang Wang} \email{zqw@fudan.edu.cn}
\affiliation{School of Mathematical Sciences, Fudan University}

	\date{\today}
	
		\maketitle

	\section{Introduction}\label{section1}
The time delay  is  ubiquitous   in many engineering control systems. Both
  actuator delay
and sensor delay  may change not only the performance
but also cause damages to the stability of the systems.
Since  the time delay is an infinite-dimensional dynamics itself,  the
  time delay in observation and control presents great
challenge in distributed parameter control systems \cite{Logemann}.
It is well known that a small
time delay in a stabilizing boundary output feedback could destabilize the system
\cite{Datko86, Datko88,Datko93}. This  implies that
some       controllers      of   PDEs  may become     practically  not implementable in the
presence of the time delay.
It is therefore  very important to consider the time delay
in  the process of the controller designs.

\textcolor{black}{In the past a few decades, many efforts have been made   for the  stabilization and \textcolor{black}{stability analysis} of PDEs
 with time delay. The seminal result on feedback stabilization for one-dimensional wave equation was considered in \cite{Datko86}}
\begin{equation} \label{20199111658}
\left\{\begin{array}{l}
  z_{tt}(x,t)=z_{xx}(x,t), x\in(0,1), t>0, \\
   z(0,t)=0,\\
 z_x(1,t)=u(t), t>0, \\
y(t)=z_t(1,t), t>0,
\end{array}\right.
\end{equation}
 where $u(t)$ is the control input and  $y(t)$ the output. 
  When the time delay is absent, the stabilization of \eqref{20199111658} is trivial
that a propositional  feedback $u(t)=-kz_t(1,t)$ with $k>0$ can stabilize exponentially the system.
However, this feedback is not robust to the input delay. In other word,
 system
\begin{equation} \label{20199111503}
\left\{\begin{array}{l}
    z_{tt}(x,t)=z_{xx}(x,t) , x\in(0,1), t>0,
 \\z(0,t)=0,
 \\z_x(1,t)=-kz_t(1,t-\tau), k>0, t\geq0
\end{array}\right.
\end{equation}
is always unstable as indicated in \cite{Datko86}, no matter how small the time delay  $\tau$ is. 
 A typical study was first made  in  \cite{XuGQ2006delay} where stabilization for
 system \eqref{20199111658} with input delay was considered by
 regarding the time delay as
   a dynamics represented by  one-dimensional   transport equation.
 The controller in \cite{XuGQ2006delay} is
split   into two parts: $u(t)=-k\mu z_t(1, t) - k(1 - \mu)z_t(1, t-\tau)$ and
the time delay free part cannot be zero. The same approach was also used in \cite{Sage1,Sage2}.
The  backstepping method is an another powerful tool to compensate the
 time delay, which has been applied  to finite-dimensional systems with actuator and sensor delays
 in \cite{Krstic2008scl}.   In monograph
 \cite[Chapter 19]{Krsticdelaybook}, stabilization for  an anti-stable wave equation
 with input delay was  realized by the backstepping method.
 Although the time delay can be  compensated in the infinite-dimensional   actuator,
  the controller,  as a full state feedback, was  very complicated in  \cite[Chapter 19]{Krsticdelaybook}.
  When the output is suffered from a time delay, an  observer/predictor
  was proposed in  \cite{GuoXuCZDelay2012}, which is  systematical
  and has been extended recently to the abstract systems  in
   \cite{MeiZD}  and \cite{GuoMEiZDTAC}.  But in general, the
   convergence in \cite{MeiZD, GuoMEiZDTAC} is valid only for smooth initial states.
   Although the convergence of \cite{GuoXuCZDelay2012} is true for general initial states,
   the controller seems still not straightforward. When
the input delay equals even multiples of the wave propagation time, the system can be
stabilized exponentially by  direct feedback discussed in \cite{Wangdelay2011}.
In \cite{Gugat}, a specific time delay in boundary observation can be used to stabilize the wave equation. \textcolor{black}{In \cite{Gugat2}, a switching delay feedback can be used to stabilize a vibrating string.} In \cite{FSIAM2} and \cite{FengTAC},  a  non-collocated feedback was
proposed to cope with the output delay being equal to one.

In this paper,  we \textcolor{black}{consider} a new controller to stabilize system \eqref{20199111658}
with input delay.  
 By considering  the   time delay  dynamic as a first order transport equation,
 the problem is converted into boundary control of a cascaded PDE system. By virtue of
boundary stabilization for first order hyperbolic systems,  we are able to design the
feedback control.
  The idea was inspired by \cite[Lemma 3.5]{WangZQ2013}  where
   the following system of transport equation
\begin{equation} \label{20189202001}
\left\{\begin{array}{l}
  \tau w_t(x,t)+w_x(x,t)=0,\\
 w(0,t)=k w(1,t) +  {(k-1) \mu}  \int_0^1w(x,t){\rm d}x
\end{array}\right.
\end{equation}
is exponentially stable
in the state space $L^2(0,1)$ if and only if $ |k |<1$ and $\mu>-1$.
 This     inspires us  that the boundary feedback   $u(t)= -c_1w(1,t) -c_2\int_0^1w(x,t){\rm d}x$ with $c_1,c_2 \in\mathbb{R}$
  may compensate   the time delay.

{\color{black}
Enlightened by \cite{Wangdelay2011}, we use the semigroup approach to explain the well-posedness of the system and the Riesz basis approach to get the dynamical behavior of the system in terms of vibrating frequencies. 
The major contribution of our paper is to develop a new method for spectral analysis.
 We derive sufficient and necessary conditions for the feedback gain and time delay which guarantee the exponential stability of the closed-loop system.  Comparing with similar conditions developed in \cite{Wangdelay2011},
 we   get 
 the \textcolor{black}{explicit term} of the stability region of $c$ for different values of $\tau$ \textcolor{black}{and from this}
  we \textcolor{black}{easily} obtain the shrink of the stability region as $\tau$ increases.
 Another main advantage of the proposed method lies in the 
  investigation of   the robustness to a small perturbation in time delay in high frequencies. 
Actually, we prove that any small perturbation of $\epsilon\in\mathbb{R}$ in time delay will excite a high frequency mode (i.e., a mode with frequency on the order of $\mathcal{O}(\dfrac{1}{|\epsilon|})$ as $\epsilon\to 0$).
  In this way, we
give an intrisic mathematical  interpretation of the destabilizing effect of arbitrarily small time delays and hence verify the judgement (or conjecture) given in \cite[Page 5, remark]{Datko86}.  This gives a mathematical explanation why numerical experiments usually do not demonstrate the non-robustness when a small perturbation is added to the time delay (See Remark \ref{robustremark}).}

We proceed as follows. 
  In Section \ref{section2}, we present the model formulation, well-posedness of the closed-loop system and fundamental spectral analysis.  In Section \ref{section3}, we derive sufficient and necessary conditions for the feedback gain and time delay which guarantee the exponential stability of the closed-loop system.  In Section \ref{section4.5}, we investigate the robustness to a small perturbation in time delay in high frequencies.  In Section \ref{section5}, we present some numerical simulations for illustration.

	\textbf{Notations}. In this paper, $\mathbb{R}\ ( \mathbb{R}_{+}, \mathbb{R}_{-})$ denote the set of all real (positive, negative) numbers, respectively. $\mathbb{C}_+\ (\mathbb{C}_-,\mathbb{C}_0)$ denote the set of complex numbers with positive (negative, zero) real parts, respectively. $\mathbb{Z}\ (\mathbb{N}^{*},\mathbb{N})$ denote the set of integers (positive integers, non-negative integers), respectively. The imaginary unit is denoted by ${\rm i}$, where $\rm i=\sqrt{-1}$. ${\rm det}(\cdot)$ denotes the determinant of a matrix. For $\lambda\in\mathbb{C}$, ${\rm Re} \lambda$, ${\rm Im} \lambda$, ${\rm arg}\lambda$ and $|\lambda|$ denote the real part, imaginary part, principle value of argument and the norm of $\lambda$, respectively. For $x\in\mathbb{R}$, ${\rm Sgn}(x)$ denotes the sign of $x$, which indicates ${\rm Sgn}(x)=1,0$ and $-1$ for $x>0$, $x=0$ and $x<0$, respectively.  For an operator $\mathscr{A}$, $D(\mathscr{A})$, $\rho(\mathscr{A})$, $\sigma(\mathscr{A})$ and $\sigma_p(\mathscr{A})$ denote the domain, regular point set, spectral point set and eigenvalue set of the operator $\mathscr{A}$, respectively. $\mathbb{D}$ denotes the unit circle $\big\{z\big||z|< 1\big\}$ in the complex plane, while $\overline{\mathbb{D}}$ denotes its closure $\big\{z\big||z|\leq 1\big\}$.

\section{Preliminaries} \label{section2}
\subsection{Model formulation}
We consider the
 stabilization of  the following wave equation with time delay in the control:
\begin{equation} \label{2018923937}
\left\{\begin{array}{l}
  z_{tt}(x,t)=z_{xx}(x,t),\\
 z(0,t)=0,\\
 z_x(1,t)=u(t-\tau),\\
y(t)=z_t(1,t),{t\geq0,}
\end{array}\right.
\end{equation}
 where $u(t)$ is the control input, $y(t)$ is the output,  and $\tau$ is the
time delay. For notational simplicity, we omit in equations hereafter the
obvious domains for both time $t$ and spatial variable $x$ when
there is no confusion.
Set $w(x,t)=u(t-\tau x)$. Then, the time delay system
\eqref{2018923937} is written as
 \begin{equation} \label{2018923934}
\left\{\begin{array}{l}
  z_{tt}(x,t)=z_{xx}(x,t),\\
 z(0,t)=0,\\
 z_x(1,t)=w(1,t),\\
 \tau w_t(x,t)+w_x(x,t)=0,\\
 w(0,t)=u(t),\\
y(t)=z_t(1,t),
\end{array}\right.
\end{equation}
which is a  {cascaded}  PDE system without explicitly the time delay.
 Inspired from \eqref{20189202001}, a  feedback control  is  designed as
\begin{equation} \label{2018923949}
\left.\begin{array}{l}
  u(t)=- c_1 w(1,t)-c_2   z_t(1,t) ,
\end{array}\right.
\end{equation}
where $c_1, c_2\in\mathbb{R}$ are   tuning parameters.  Since we only consider the
output feedback, the integral term $\int_0^1w(x,t){\rm d}x$  in  \eqref{20189202001}  is ignored.
The first term of \eqref{2018923949} is damper for the transport equation
and the second term is a direct proportional feedback for wave equation without time delay. When $c_1=0$, the system becomes \eqref{20199111503}.
Under the feedback  \eqref{2018923949}, the closed-loop of system \eqref{2018923934} reads
\begin{equation} \label{20189202158}
\left\{\begin{array}{l}
  z_{tt}(x,t)=z_{xx}(x,t),\\
 z(0,t)=0,\\
 z_x(1,t)=w(1,t),\\
 \tau w_t(x,t)+w_x(x,t)=0,\\
 w(0,t)= kz_t(1,t) .
\end{array}\right.
\end{equation}

	\subsection{Well-posedness of system \eqref{20189202158}}

We consider system \eqref{20189202158} in the state space
\begin{equation} \label{20199121519}
 \mathcal{X}=H^1_L(0,1)\times (L^2(0,1))^2,
 \end{equation}
  where
$H^1_L(0,1)=\{f\in H^1 (0,1)\ |\ f(0)=0\}$. The inner product in $ \mathcal{X}$ is defined  by
\begin{equation} \label{20189241648}
\langle (f_1,g_1,h_1),(f_2,g_2,h_2)\rangle_{ \mathcal{X}}=\int_0^1[ f_1'(x)\overline{ f_2'(x)}
+ g_1(x)\overline{g_2(x)}+  h_1(x)\overline{h_2(x)}]{\rm d}x
\end{equation}
for $   (f_i,g_i,h_i)\in \mathcal{X},\ i=1,2 $.
System \eqref{20189202158} can be written as an evolutionary equation in $ \mathcal{X}$:
\begin{equation} \label{20189241451}
\dfrac{\rm d}{{\rm d}t}(z(\cdot,t),z_t(\cdot,t),w (\cdot,t))= \mathscr{A} (z(\cdot,t),z_t(\cdot,t),w (\cdot,t)),
\end{equation}
where the operator $ \mathscr{A}: D( \mathscr{A})\subset  \mathcal{X}\to \mathcal{X}$ is defined by
\begin{equation} \label{20189241453}
\left\{\begin{array}{l}
 \mathscr{A}(f,g,h)=(g,f'',-\tau^{-1}h' ),\forall (f,g,h)\in D(\mathscr{A}),
\\D(\mathscr{A})=\left\{( f,g,h) \ |\ f\in H^2(0,1),g,h\in H^1(0,1),f(0)=g(0)=0,\right.
\\ \hspace{2cm}\left. f'(1)=h(1),\ h(0)=-c_1 h(1)-c_2 g(1)\right\}.
\end{array}\right.
\end{equation}

\begin{lemma}\label{Lemma20232101235}
Let  the operator  $\mathscr{A} $ be defined by  \eqref{20189241453}. Then $\mathscr{A}^{-1}$ exists and is compact. Hence, $\sigma(\mathscr{A})$ consists of isolated eigenvalues multiplicity only.
\end{lemma}
Proof
For any $(\hat{f},\hat{g},\hat{h})\in\mathcal{X}$, we solve the equation $\mathscr{A}(f,g,h)=(\hat{f},\hat{g},\hat{h})$ to get
 \begin{equation} \label{20189261005-2}
\left\{\begin{array}{l}
  g=\hat{f},
\\ f(x)=f'(1)x-\int_0^x\int_{\alpha}^1\hat{g}(s){\rm d}s,\ f'(1)=-\frac{1}{1+c}
\left( c \hat{f}(1)+ \tau\int_0^1\hat{h}(s){\rm d}s\right),
\\ h(x)=h(0)-\tau\int_0^x\hat{h}(s){\rm d}s,\ h(0)=\frac{c}{1+c }
\left(  \tau\int_0^1\hat{h}(s){\rm d}s-  \hat{f}(1)\right),
\end{array}\right.
\end{equation}
which implies, by the Sobolev trace-embedding,  that $\mathscr{A}^{-1}$
exists and is compact in $\mathcal{X}$. So
$\sigma(\mathscr{A})$ consists of isolated eigenvalues of finite algebraic multiplicity.
~~

\begin{theorem}\label{wellposedness}
Suppose that  $c_1, c_2\in \mathbb{R} $ and   $  \tau  >0$.
  Then,  the operator  $\mathscr{A} $  defined by  \eqref{20189241453} generates a
$C_0$-group ${\rm e}^{\mathscr{A}t}$ on $\mathcal{X}$.
\end{theorem}
Proof
Inspired by  \cite[Theorem 2.2]{Wangdelay2011},
we first introduce  a new inner product
 \begin{small}\begin{equation*} \label{20189241453*0}
\left.\begin{array}{l}
     \langle(f_1,g_1,h_1),(f_2,g_2,h_2)\rangle_1
\triangleq \int_0^1 {\rm e}^{\alpha x}(f_1'-g_1)\overline{ (f_2'-g_2)}{\rm d}x
+ \int_0^1 {\rm e}^{\beta x}(f_1'+g_1)\overline{( f_2'+g_2)}{\rm d}x \\ 
     \hspace{46mm} +\tau \int_0^1 {\rm e}^{\gamma x} h_1(x)\overline{h_2(x)}{\rm d}x ,\quad \
\forall\ (f_i,g_i,h_i)\in\mathcal{X},\ i=1,2 ,
\end{array}\right.
\end{equation*}\end{small}
where
\begin{equation} \label{201998920}
  \beta\leq 0,\ \ {\rm e}^{\alpha}>1+c_2^2,\ \ {\rm e}^{\gamma}>\frac{\left({\rm e}^{\alpha}+{\rm e}^{\beta}+|c_1c_2|\right)^2}{{\rm e}^{\alpha}-1-c_2^2}.
  \end{equation}
For $( {f}, {g}, {h})\in D(\mathscr{A})$, a simple computation shows that
 \begin{small}
 \begin{align}
 &{\rm Re}\langle \mathscr{A}(f,g,h),  (f,g,h) \rangle_{1}\\ \notag
 =&{\rm Re} \langle (g,f'',-\tau^{-1}h' ),  (f,g,h) \rangle_{1}
  \\=& -\frac12{\rm e}^{\alpha x}|g-f'|^2\Big{|}_0^1+\frac{\alpha}{2}\int_0^1
 {\rm e}^{\alpha x}|g-f'|^2{\rm d}x +
   \frac12{\rm e}^{\beta x}|g+f'|^2\Big{|}_0^1-\frac{\beta}{2}\int_0^1
 {\rm e}^{\beta x}|g+f'|^2{\rm d}x \notag
 \\& -\frac12{\rm e}^{\gamma x}|h|^2\Big{|}_0^1+\frac{\gamma}{2}\int_0^1
 {\rm e}^{\gamma x}|h |^2{\rm d}x \notag
  \\=&-\frac12{\rm e}^{\alpha}|g(1)-h(1)|^2+ \frac12{\rm e}^{\beta}|g(1)+h(1)|^2
 -\frac12{\rm e}^{\gamma}|h(1)|^2+ \frac12 | h(0)|^2 \notag
   \\ &+\frac{\alpha}{2}\int_0^1
 {\rm e}^{\alpha x}|g-f'|^2{\rm d}x  -\frac{\beta}{2}\int_0^1
 {\rm e}^{\beta x}|g+f'|^2{\rm d}x+\frac{\gamma}{2}\int_0^1
 {\rm e}^{\gamma x}|h |^2{\rm d}x \notag
   \\=&-\frac12{\rm e}^{\alpha}|g(1)-h(1)|^2+ \frac12{\rm e}^{\beta}|g(1)+h(1)|^2+\frac12|c_1h(1)+c_2g(1)|^2
 -\frac12{\rm e}^{\gamma}|h(1)|^2 \notag
   \\& +\frac{\alpha}{2}\int_0^1
 {\rm e}^{\alpha x}|g-f'|^2{\rm d}x  -\frac{\beta}{2}\int_0^1
 {\rm e}^{\beta x}|g+f'|^2{\rm d}x+\frac{\gamma}{2}\int_0^1
 {\rm e}^{\gamma x}|h |^2{\rm d}x. \notag
      \end{align}
\end{small}
  Thus
  \begin{small}
 \begin{align} \label{20189261005-1}
 & {\rm Re}\langle \mathscr{A}(f,g,h),  (f,g,h) \rangle_{1}
  \\ \leq & -\frac12\left({\rm e}^{\alpha}-{\rm e}^{\beta}-c_2^2\right)|g(1)|^2-
 \frac12\left({\rm e}^{\alpha}+{\rm e}^{\gamma}-{\rm e}^{\beta}-c_1^2\right)|h(1)|^2 \notag
   \\&+\left({\rm e}^{\alpha}
 +{\rm e}^{\beta}+|c_1c_2|\right)|h(1)g(1)|
    +\frac{\alpha}{2}\int_0^1{\rm e}^{\alpha x}|g-f'|^2{\rm d}x-\frac{\beta}{2}\int_0^1
 {\rm e}^{\beta x}|g+f'|^2{\rm d}x \notag
  \\&+\frac{\gamma}{2}\int_0^1
 {\rm e}^{\gamma x}|h |^2{\rm d}x. \notag
    \end{align}
\end{small}
By Young's inequality, for any $\delta>0$,
 \begin{equation} \label{2019911824}
|h(1)g(1)|\leq \frac{\delta}{2}|g(1)|^2+\frac{1}{2\delta}|h(1)|^2.
\end{equation}
Combine \eqref{20189261005-1} and \eqref{2019911824} to obtain
\begin{equation} \label{20189261005}
 \begin{array}{l}
 {\rm Re}\langle \mathscr{A}(f,g,h),  (f,g,h) \rangle_{1}
 \leq -a_1|g(1)|^2-a_2|h(1)|^2
      +\frac{\alpha}{2}\int_0^1\\
 {\rm e}^{\alpha x}|g-f'|^2{\rm d}x  -\frac{\beta}{2}\int_0^1
 {\rm e}^{\beta x}|g+f'|^2{\rm d}x+\frac{\gamma}{2}\int_0^1
 {\rm e}^{\gamma x}|h |^2{\rm d}x,
    \end{array}
\end{equation}
where
\begin{equation} \label{20189911941}
\left\{\begin{array}{l}
a_1\triangleq \frac12\left[{\rm e}^{\alpha}-{\rm e}^{\beta}-c_2^2-\delta\left({\rm e}^{\alpha}+{\rm e}^{\beta}+|c_1c_2|\right)\right],\\
a_2\triangleq  \frac12\left[{\rm e}^{\alpha}+{\rm e}^{\gamma}-{\rm e}^{\beta}-c_1^2-\frac{{\rm e}^{\alpha}+{\rm e}^{\beta}+|c_1c_2|}{\delta}\right].
\end{array}\right.
\end{equation}
Owing to  \eqref{201998920},  we  can choose  $\delta$ small enough such that
\begin{equation} \label{2019926811}
0<\delta \left({\rm e}^{\alpha}+{\rm e}^{\beta}+|c_1c_2|\right) <{\rm e}^{\alpha}-1-c_2^2.
\end{equation}
By \eqref{201998920}, it follows that
\begin{equation} \label{2019926812}
 \frac{{\rm e}^{\alpha}+{\rm e}^{\beta}+|c_1c_2|}{\delta}<{\rm e}^{\gamma}+{\rm e}^{\alpha}-1-c_1^2
 \leq {\rm e}^{\gamma}+{\rm e}^{\alpha}-{\rm e}^{\beta}-c_1^2.
\end{equation}
By \eqref{20189911941}, \eqref{2019926811} and \eqref{2019926812},
\begin{equation} \label{20189911957}
a_1>0\ \ \mbox{ and}\ \  a_2>0.
\end{equation}
Hence, it follows from \eqref{20189261005}  and   \eqref{20189911957}  that
  there is an $M  > 0 $ such that
   $\forall\  (f,g,h) \in D(\mathscr{A})$,
  \begin{equation} \label{20199111019}
 {\rm Re}\langle \mathscr{A}(f,g,h),  (f,g,h) \rangle_{1}
 \leq  -a_1|g(1)|^2-a_2|h(1)|^2 +
 M\langle (f,g,h),  (f,g,h) \rangle_{1},
\end{equation}
This implies that $\mathscr{A} - M $ is dissipative. By \eqref{20189261005}, $\mathscr{A}$ is a discrete operator
 (i.e., $\mathscr{A}^{-1}$ is
compact), so there is a sequence $M_n\to\infty$ such that $M_n\in \rho(\mathscr{A})$, the resolvent set of
$\mathscr{A}$. We may assume without loss of generality that $M \in \rho(\mathscr{A})$. By the Lumer-Phillips theorem, $\mathscr{A}- M$ generates a $C_0$-semigroup of contractions
${\rm e}^{(\mathscr{A}-M)t}$ on  $\mathcal{X}$ (\cite[Theorem 4.3, p.14]{Pazy}). The bounded perturbation theorem of $C_0$-semigroups
ensures that $\mathscr{A}$ generates a $C_0$-semigroup ${\rm e}^{\mathscr{A}t}$ on $\mathcal{X}$
(\cite[Theorem 1.1, p.76]{Pazy}).
Similarly, we apply the same argument  to get that
$-\mathscr{A}$ also generates a $C_0$-semigroup in $\mathcal{X}$ (see, e.g., \cite{Wangdelay2011}).
Therefore, $\mathscr{A}$  actually generates a $C_0$-group
on $\mathcal{X}$. This completes the proof.
~~

\subsection{Spectral analysis}
Let us now consider the eigenvalue problem of $\mathscr{A}$.

Let us now consider the eigenvalue problem of $\mathscr{A}$, $\mathscr{A}(f,g,h)=\lambda(f,g,h)$, then
  $g=\lambda f$ and
\begin{equation} \label{20189231117}
\left\{\begin{array}{l}
  f''=\lambda^2f, h'=-\tau\lambda h,\\
f(0)=0, f'(1)=h(1),\\
 h(0)=-c_1 h(1)-   c_2 \lambda f(1)  .
\end{array}\right.
\end{equation}

When $\lambda=0$, the solution of \eqref{20189231117} is found to be
\begin{equation}
\left\{\begin{array}{l}
f(x)=bx,
\\h(x)=b,
 \end{array}\right.
 \end{equation}
where the constant $b$ satisfies that $b=-c_1b$. Thus, the equation has a nontrivial solution if and only if $c_1=-1$. When $c_1=-1$, the corresponding eigenfunction
$(f_0, g_0, h_0)$ is given by
$$
f_0(x)=x,~~  g_0(x)=1, ~~  h_0(x)=1.
$$

When $\lambda\neq 0$, the solution of \eqref{20189231117} is found to be
\begin{equation} \label{20189231123}
\left\{\begin{array}{l}
  f(x)=a{\rm e}^{\lambda x}-a{\rm e}^{-\lambda x}=2a\sinh \lambda x,\\
  h(x)=b{\rm e}^{-\tau\lambda x},
 \end{array}\right.
\end{equation}
where $a$ and $b$ are constants that satisfy
\begin{equation} \label{20189231125}
\left\{\begin{array}{l}
  c_2\lambda\left({\rm e}^{\lambda}-{\rm e}^{-\lambda}\right) a
+\left(1+c_1 {\rm e}^{-\tau\lambda} \right) b=0,\\
     \lambda\left({\rm e}^{\lambda}+{\rm e}^{-\lambda}\right) a
  -{\rm e}^{-\tau\lambda} b=0.
 \end{array}\right.
\end{equation}
 The characteristic determinant of \eqref{20189231125}
is
\begin{equation} \label{20189231135}
 \begin{array}{l}
 \Delta(\lambda) =\left|\begin{array}{cc}
    c_2 \lambda\left({\rm e}^{\lambda}-{\rm e}^{-\lambda}\right)&
  1+c_1{\rm e}^{-\tau\lambda} \\
    \lambda\left({\rm e}^{\lambda}+{\rm e}^{-\lambda}\right) &
  -{\rm e}^{-\tau\lambda}
 \end{array}\right|
  =-\lambda\cosh \lambda(1+c_1{\rm e}^{-\lambda\tau})-c_2\lambda\sinh \lambda {\rm e}^{-\lambda\tau}. \end{array}
\end{equation}
Thus, the equation has a nontrival solution  if and only if
 \begin{equation} \label{20189231232}
\tilde{\Delta}(\lambda)\triangleq-\cosh \lambda(1+c_1{\rm e}^{-\lambda\tau})-c_2\sinh \lambda {\rm e}^{-\lambda\tau}=0.
\end{equation}
Furthermore, we observe that when $c_1=-1$, $\lambda=0$ is also a root for $\tilde{\Delta}(\lambda)=0$, then we conclude that
$$
\sigma({\mathscr{A}})=\sigma_p(\mathscr{A})=\{\lambda\in\mathbb{C}|\tilde{\Delta}(\lambda)=0\}.
$$
Therefore, each $\lambda\in\sigma(\mathscr{A})\backslash\{0\}$ is geometrically simple, and the corresponding eigenfunction
$(f_\lambda, g_\lambda, h_\lambda)$ is given by
$$
f_\lambda(x)={\rm e}^{-\tau\lambda}\sinh \lambda x,~~  g_\lambda(x)=\lambda{\rm e}^{-\tau\lambda}\sinh \lambda x, ~~  h_\lambda(x)=\lambda\cosh \lambda {\rm e}^{-\tau\lambda x}.
$$
For $\lambda=0, c_1=-1$, the corresponding eigenfunction
$(f_0, g_0, h_0)$ is given by
$$
f_0(x)=x,~~  g_0(x)=1, ~~  h_0(x)=1.
$$

For any $\lambda\in\rho(\mathscr{A})$, we have following lemma on the expression of the resolvent operator.
\begin{lemma}\label{Lemma20199101635}
Let  the operator  $\mathscr{A} $ be defined by  \eqref{20189241453}.
 Then,
for any $\lambda\in \rho(\mathscr{A})$ and
$Y=(f_1,g_1,h_1)\in {\mathcal X}$, $X=R(\lambda,\mathscr{A})Y$,
where $X=(f,g,h)\in D(\mathscr{A})$ is given by
\begin{equation}\label{wxh2019992018}\left\{\begin{array}{l}
f(x,\lambda)=laystyle\frac{(1+c_1{\rm e}^{-\tau\lambda})F_1+F_2{\rm e}^{-\tau\lambda}}
{ \Delta(\lambda)}\sinh\lambda x+F_0(x,\lambda),\\
g(x,\lambda)=\lambda f(x,\lambda)-f_1(x),\\
h(x,\lambda)=laystyle\frac{F_2\lambda \cosh\lambda-F_1c_2\lambda\sinh\lambda}{ \Delta(\lambda)}{\rm e}^{-\tau\lambda x}+H_0(x,\lambda),
\end{array}\right.\end{equation}
and
\begin{equation}\label{wxh2019992059}\left\{\begin{array}{l}
  \Delta(\lambda)=\lambda\cosh \lambda(1+c_1{\rm e}^{-\lambda\tau})+c_2\lambda\sinh \lambda {\rm e}^{-\lambda\tau},\\
F_0(x,\lambda)=laystyle
\int_{0}^{x}\lambda^{-1}\sinh\lambda(x-s)[-\lambda f_1(s)-g_1(s)]{\rm d}s,\\
H_0(x,\lambda)=laystyle\int_{0}^{x}{\rm e}^{-\tau\lambda(x-s)}\tau h_1(s){\rm d}s,\\
F_1=laystyle\int_{0}^{1}\cosh\lambda(1-s)[\lambda f_1(s)+g_1(s)]{\rm d}s+
\int_{0}^{1}{\rm e}^{-\tau\lambda(1-s)}\tau h_1(s){\rm d}s,\\
F_2=laystyle\int_{0}^{1} c\sinh\lambda(1-s)[\lambda f_1(s)+g_1(s)]{\rm d}s-
\int_{0}^{1}c{\rm e}^{-\tau\lambda(1-s)}\tau h_1(s){\rm d}s.
\end{array}\right.\end{equation}
\end{lemma}
Proof:
For any $Y=(f_1,g_1,h_1)\in{\mathcal X}$ and  $\lambda\in\rho(\mathscr{A})$,
let
\begin{equation}\label{wxh2019992020}
X=R(\lambda,\mathscr{A})Y, \ X=(f,g,h)\in D(\mathscr{A}).
\end{equation}
Then,
\begin{equation}\label{wxh2019992022}
(\lambda I-\mathscr{A})X=(\lambda f-g,\lambda g-f^{\prime\prime},\lambda h+\tau^{-1}h')=(f_1,g_1,h_1).
\end{equation}
Hence,  $g=\lambda f-f_1$ with $f,h$ satisfying
\begin{equation}\label{wxh2019992023}
\left\{\begin{array}{l}
f^{\prime\prime}(x)-\lambda^2 f(x)=-g_1-\lambda f_1,\\
f'(1)=h(1),\\
h^{\prime}(x)+\tau\lambda h=\tau h_1,\\
h(0)=-c_1h(1)-c_2\lambda f(1),
\end{array}\right.
\end{equation}
which gives
\begin{equation}\label{wxh2019002045}
\left\{\begin{array}{l}
f(x,\lambda)=a\sinh\lambda x+F_0(x,\lambda),\\
h(x,\lambda)=b {\rm e}^{-\tau\lambda x}+H_0(x,\lambda),
\end{array}\right.
\end{equation}
where $F_0(x,\lambda)$ and $H_0(x,\lambda)$ are given by
\eqref{wxh2019992059}.
By the boundary conditions of
\eqref{wxh2019992023}, we have
\begin{equation}\label{wxh2019992208}\left\{\begin{array}{l}
a\lambda\cosh\lambda-b{\rm e}^{-\tau\lambda}=F_1,\\
ac_2\lambda\sinh\lambda+b(1+c_1{\rm e}^{-\tau\lambda})=F_2,
\end{array}\right.a,b\in\mathbb{R},
\end{equation}
where $F_1$ and $F_2$ are given by
\eqref{wxh2019992059}. Now, we determine the constants $a$ and $b$.
Since   $\lambda\in\rho(\mathscr{A})$ and $\Delta(\lambda) $ happens to be the characteristic determinant
of \eqref{wxh2019992208}, it follows that $ \Delta(\lambda)\neq 0$. Therefore,
 $a$ and $b$ can be determined by solving
 equation \eqref{wxh2019992208}.
Moreover, the solution $X$ of \eqref{wxh2019992020} can
be written in \eqref{wxh2019992018}.

$$
\Delta(\lambda)={\rm e}^{-\lambda}+k{\rm e}^{-(1+\tau)\lambda}+{\rm e}^{\lambda}-k{\rm e}^{(1-\tau)\lambda}
$$

Enlightened by \cite[Proposition 3.3]{Wangdelay2011}, we characterize the spectrum of $\mathscr{A}$ as follows.
\begin{theorem}\label{20232101439}
Let $\mathscr{A} $  defined by  \eqref{20189241453} and $\tilde{\Delta}(\lambda)$ defined by \eqref{20189231232}. The following assertions hold for the spectrum of $\mathscr{A}$:
\\(i) There is an $M>0$ such that for all $\lambda\in\rho(\mathscr{A})$, $|{\rm Re}\lambda|<M$; that is, all the eigenvalues of $\mathscr{A}$ lies in some vertical strip paraller to the imaginary axis in the complex plane.
\\(ii)The multiplicity of each root of $\tilde{\Delta}(\lambda)=0$ is at most two.
\\(iii) If $\tau$ is rational, then the eigenvalue of $\mathscr{A}$ are located on finitely many lines parallel to the imaginary axis.
\\(iv)If $\tau$ is irrational, then all roots of $\tilde{\Delta}(\lambda)=0$ are simple.
\\(v)The eigenvalues of $\mathscr{A}$ are separated, that is
$$
\inf_{\lambda_m,\lambda_n\in\sigma(\mathscr{A}),\lambda_m\neq\lambda_n}|\lambda_m-\lambda_n|>0.
$$
\\(vi) The algebric multiplicity of each eigenvalue of $\mathscr{A}$ is at most two.
\end{theorem}
The proof of Theorem \ref{20232101439} is similar to  that of \cite[Proposition 3.3]{Wangdelay2011}, so we omit it.

\begin{lemma}\label{Lm20199101733}
 Let  $\mathscr{A} $       be given by \eqref{20189241453}. Then, the root subspace of $\mathscr{A} $  is complete
in $\mathcal{X}$, that is, ${\rm Sp}( \mathscr{A}   ) = \mathcal{X}$, where ${\rm Sp}( \mathscr{A}   ) $ denotes the root subspace of  $ \mathscr{A}  $  spanned by the
generalized eigenfunctions of $ \mathscr{A}  $.
 \end{lemma}
Proof:
From Lemma \ref{Lemma20199101635}, $X = R(\lambda,\mathscr{A} )Y$ can be further represented as
$$
 X= R(\lambda,\mathscr{A} )Y=\frac{G(\lambda,Y )}{\Delta(\lambda)},
$$
where $G(\lambda,Y )$ is an  $\mathcal{X}$-valued entire function with order less than or equal to 1, and
by \eqref{20189231135}, $\Delta(\lambda)$ is a scalar entire function of order 1. Since  from  Theorem \ref{wellposedness}, $\mathscr{A}$ generates a $C_0$-group on $\mathcal{X}$,
$|R(\lambda,\mathscr{A})|$ is uniformly bounded as ${\rm Re}\lambda\to \pm\infty$.
 By \cite[Theorem 4.1]{Wangdelay2011} or \cite[Theorem 4]{XuGQ2003}, ${\rm Sp}( \mathscr{A}   ) = \mathcal{X}$. This completes the
  proof of the lemma.

\begin{lemma}\label{Lm20199101844}
Let  the operator  $\mathscr{A} $ be defined by  \eqref{20189241453}.
Then,  the  spectrum determined growth condition holds for $\mathscr{A} $: $s(\mathscr{A} ) = \omega(\mathscr{A} )$,
where $s(\mathscr{A} )$ and  $\omega(\mathscr{A} )$ are the spectral bound of $\mathscr{A}$ and the
the growth order of ${\rm e}^{\mathscr{A}t}$, respectively.
\end{lemma}
Proof
  By \eqref{20189261005}, $\mathscr{A}$ is a discrete operator. Suppose that
  $\{\lambda_n\}_{n=1}^{\infty}$ is the set of the eigenvalues of $\mathscr{A}$.
 Then,  $\tilde{\Delta}(\lambda_n)=0$, where the function $\tilde{\Delta}$ is defined by \eqref{20189231232}.
  It is evident that $\tilde{\Delta}(\lambda) $ is   an entire function of
exponential type. Moreover, by Theorem \ref{20232101439}   and the fact that
$|\tilde{\Delta}(\lambda)| \to\infty$ as
${\rm Re}\lambda \to \pm\infty$, it follows that $\tilde{\Delta}(\lambda) $ is a sine-type function.
  By \cite[Theorem 1]{XuGQ2003} or \cite{PSulian1979},
   $\{{\rm e}^{\lambda_nt}\}_{n=1}^{\infty}$ forms a Riesz basis for $L^2(0, T)$ for some $T > 0$.
Combining Theorem \ref{20232101439}-(v), Lemma \ref{Lm20199101733}  and
\cite[Theorem 4.3]{Wangdelay2011}, the spectrum-determined growth condition holds for $\mathscr{A} $.

A necessary and sufficient condition that the polynomial 
$$
F(\lambda)=a_m\lambda^m+a_{m-1}\lambda^{m-1}+\cdots+a_1\lambda+a_0, a_m>0,
$$
with real coefficient have all of its roots inside the unit circle is given by
$$
F(1)>0,~~~~~ (-1)^{m}F(-1)>0,
$$
and the $(m-1)\times (m-1)$ Jury matrices
$$
\Delta^{\pm}_{m-1}=\left(\begin{array}{lllll}
   a_m &0&0&\cdots&0
      \\  a_{m-1} & a_m & 0 &\cdots &0
      \\ a_{m-2} & a_{m-1} & a_m &\cdots & 0
      \\ \vdots & \vdots & \vdots & \vdots & \vdots
      \\ a_2 & a_3 & a_4 &\cdots &a_m
\end{array}\right)\pm \left(\begin{array}{lllll}
   0 &0&\cdots&0&a_0
      \\  0 &0 & \cdots &a_0 & a_1
      \\ \vdots & \vdots& \vdots &\vdots& \vdots
      \\ 0& a_0 & \cdots & a_{m-4} & a_{m-3}
      \\ a_0 & a_1 & \cdots &a_{m-3} &a_{m-2}
\end{array}\right)$$
are both positive innerwise; that is, the determinants of all the inners of $\Delta^{\pm}_{m-1}$ are positive. Here, the inners of a square matrix are the matrix itself and all the
matrices obtained by omitting successively the first and last rows and the first and last
columns.

\section{Exponential stability of system \eqref{20189202158}}\label{section3}
In this section, we discuss the exponential stability of the system \eqref{20189202158}. By Lemma \ref{Lm20199101844}, we need only to check whether all the eigenvalues of $\mathscr{A}$ is located in $\mathbb{C}_{-}$.  By $\tilde{\Delta}(\lambda)=0$, we obtain
\begin{equation}\label{cha}
\tilde{\Delta}(\lambda)\triangleq-\cosh \lambda(1+c_1{\rm e}^{-\lambda\tau})-c_2\sinh \lambda {\rm e}^{-\lambda\tau}=0.
\end{equation}
For the sake of simplicity, we consider the case $c_1=c_2=c\in\mathbb{R}$ in this section. The stability region on the whole $c_1-c_2$ plane will be investigated in our future work.
Thus, Eq. \eqref{cha} becomes
\begin{equation} \label{202302042155}
{\rm e}^{2\lambda}+2c {\rm e}^{(2-\tau)\lambda}+ 1=0.
\end{equation}

\subsection{$\tau>0$ is rational}
First, we talk about the situation that $\tau>0$ is rational.
Eq. \eqref{202302042155} can be written as
$
c=-\dfrac{1}{2}[{\rm e}^{\tau\lambda}+{\rm e}^{(\tau-2)\lambda}].
$

If $\tau=1$, suppose that $z={\rm e}^{{\lambda}}$, the Eq. \eqref{202302042155} becomes $z^2+2cz+1=0$.  No matter what the value of $c$ is, this equation
 has at least one root satisfying  $|z|\geq 1$, which indicates that Eq. \eqref{202302042155} has at least one root located in $\mathbb{C}_0\cup\mathbb{C}_{+}$.
$$
\tau=\dfrac{m}{n}$$

$$z={\rm e}^{-{\lambda}/{n}}$$
 
$$
z^n+kz^{m+n}+z^{-n}-kz^{m-n}=0
$$
no roots are inside the unit circle.

The characteristic equation \eqref{202302042155} has no complex root located in $\mathbb{C}_0\cup\mathbb{C}_{+}$ is equivalent to that the equation $f(z)=c$ has no root which is located in the unit circle $\overline{\mathbb{D}}$.

Denote by the set
$$
\mathscr{E}=\left\{k\in\mathbb{R} \Bigg| \mbox{there exists}~z\in\mathbb{C}, |z|=1~~\mbox{such that}~ f(z,k)=0 \right\}.
$$
For $k\notin\mathscr{E}_{m,n}$, we denote by a function $N(k)$, which is the number of root, counted by multiplicity, of the equation $f(z,k)=0$ located in the unit circle $\mathbb{D}$.

In order to study the stability region, we  find the region for $c$ in which $N(c)=0$. We  give an outline of our following dissertation.  $N(c)$ is the number of root branches for different values of $c$.  Each time when $c$ moves across some critical value $c\in\mathscr{E}$, a root branch disappear or emerge at the boundary of unit circle. We use implicit theorem to investigate whether a root branch disappears or emerges. When $\tau<1$, we move $c$ from $\pm\infty$ to $0$. When $c$ is at $\pm\infty$, we prove that there are $m$ branches of root by argument principle (Lemma \ref{lemma6}). Each time when $c$ comes arcoss a critical value $c$, we prove that a new branch emerges  from the boundary of the unit circle. As a result, we prove that there is no stability region. When $\tau>1$, we move $c$ from $0$ to $\pm\infty$. We prove that each time when $c$ comes arcoss a critical value $c\neq 0$, a new branch emeges from the boundary of the unit circle. Thus, we only need to investigate the root branch near $c=0$. When $c=0$, there are several roots on the unit circle. We prove that only when $\tau$ is an even number, there is a region for $c$ such that all the roots leave the unit circle and therefore there is a stability region for $c$.

According to Lemma \ref{lemma6}, we get that for an interval $[c_1,c_2]$, if $[c_1,c_2]\cap \mathscr{E}=\emptyset$, $N(c)$ is a constant. And when $|c|>1$, $N(c)\equiv m$.

Now we turn to investigate  the relationship between $\lim_{c\to c_*-}N(c)$ and  $\lim_{c\to c_*+}N(c)$ for any $c_*\in\mathscr{E}$. We will discuss into following categories. Case A: $\tau<1$. Case B: $\tau>1$.

Case A: $\tau<1$. According to Lemma \ref{lemma2}, for each $c_*=f(z_*)\in\mathscr{E}\backslash\{0\}, |z_*|=1$, there exists an implicit function $z(c)$ such that $z(c_*)=z_*$ and $c=f(z(c))$ for each $c\in(c_*-\epsilon, c_*+\epsilon)$. Here, $\epsilon>0$ is efficiently small. Suppose that $z(c)=r(c){\rm e}^{{\rm i}\theta(c)}$, where $r(c)$ and $\theta(c)$ represent the absolute value and argument value of the function $z(c)$, respectively. Then we have
$$
{\rm Sgn}[r'(c_*)]={\rm Sgn}(n-m){\rm Sgn}(c_*)={\rm Sgn}(c_*).
$$
When $c_*>0, r'(c_*)>0$. We obtain that $|z(c)|>1$ for $c>c_*$ and $|z(c)|<1$ for $c<c_*$. This indicates that when $c$ goes from $c_*-$ to $c_*+$, a root for the equation $f(z)=c$ enters the unit circle from the outside. Thus we obtain that $\lim_{c\to c_*-}N(c)>\lim_{c\to c_*+}N(c)$ for $c>0$. Similarly, we obtain that when $c_*<0$, $\lim_{c\to c_*-}N(c)<\lim_{c\to c_*+}N(c)$. Therefore, $N(c)$ is monotonically increasing on $\mathbb{R}_{-}\backslash\mathscr{E}$ and monotonically decreasing on $\mathbb{R}_+\backslash\mathscr{E}$. Since $N(c)=m$ when $|c|>1$, $N(c)$ is not zero on $\mathbb{R}\backslash\mathscr{E}$.

Case B: $\tau>1$. Similarly, we get that $N(c)$ is monotonically decreasing on $\mathbb{R}_{-}\backslash\mathscr{E}$ and monotonically increasing on $\mathbb{R}_+\backslash\mathscr{E}$. We try to find the region for $c$ in which $N(c)=0$. Thus, we need to investigate $N(c)$ near $c=0$.

When $c=0$, $f(z)=0$ has $2n$ roots which read
$$
z_k={\rm e}^{\rm i\theta_k}, \theta_k=\dfrac{(2k+1)\pi}{2n} ,k=0,1,2,\cdots, 2n-1.
$$
There exist $2n$ different implicit functions $z_k(c)$ such that $z_k(0)=z_k={\rm e}^{\rm i \theta_k}$ with $\theta_k=\dfrac{(2k+1)\pi}{2n}$, $k=0,1,2,\cdots, 2n-1$. According to Lemma \ref{lemma2}, we have ${\rm Sgn}[r_k'(0)]={\rm Sgn}\cos[\dfrac{m(2k+1)\pi}{2n}]$. We discuss into following categories. Case I: $n>1$. Case  II: $n=1, m\geq 3$ is an odd number. Case  III: $n=1, m=4s-2, s\in\mathbb{N}^*$. Case IV: $n=1, m=4s, s\in\mathbb{N}^*$.

Case  I: $n>1$.  Lemma \ref{lemma3} shows that there exist at least two integer numbers $k,j$ such that $r_k'(0)>0$ and $r_j'(0)<0$. This implies that there exists small efficiently $\delta>0$, $|z_k(c)|<1$ for $c\in(-\delta,0)$ and  $|z_j(c)|<1$  for   $c\in(0,\delta)$. Thus $\lim_{c\to 0+}N(c)$ and $\lim_{c\to 0-}N(c)$ are both nonzero.  As a result, $N(c)$ is nonzero on $\mathbb{R}\backslash\mathscr{E}$.

Case II: $n=1$, $m\geq 3$ is an odd number. According to Lemma \ref{lemma5}, we obtain that $r_k'(0)=0$ and $r_k''(0)<0$. By Taylor expression
$
r_k(c)=r_k(0)+r_k'(0)c+\dfrac{1}{2}r_k''(0)c^2+o(c^2),
$
we obtain that $|z_k(c)|<1$ for $c$ in a neighborhood of $0$. This implies that $N(c)$ is nonzero in the neighborhood of $0$. Thus, $N(c)$ is nonzero on $\mathbb{R}\backslash\mathscr{E}$.

Case III: $n=1$, $m=4s-2, s\in\mathbb{N}^*$. A direct computation leads to that $r_k'(0)<0$ for both $k=0,1$. This indicates that $|z_k(c)|<1$ when $c\to 0+$ and $|z_k(c)|>1$ when $c\to 0-$. This further indicates that $N(c)$ is nonzero for $c\to 0+$ and $N(c)$ is zero when $c\to 0-$.  Lemma \ref{lemma4} tells us the nearest element    in $\mathscr{E}\cap\mathbb{R}_{-}$      to $0$  is  $-\sin [\dfrac{\pi}{2(m-1)}]$. Therefore, we obtain that $N(c)=0$ on the interval $(-\sin [\dfrac{\pi}{2(m-1)}],0).$

Case IV: $n=1$, $m=4s, s\in\mathbb{N}^*$. Using the same argument, we get that $N(c)=0$ on the interval $(0,\sin [\dfrac{\pi}{2(m-1)}])$.

\subsection{$\tau>0$ is irrational}\label{subsection32}
In this subsection, we prove that for any irrational $\tau>0$ and $c\in\mathbb{R}$, Eq. \eqref{202302042155} has at least one root located in $\mathbb{C}_0\cup\mathbb{C}_+$. The idea of the proof is similar to that in the last subsection. In this case, Eq. \eqref{202302042155} with respect to $\lambda$ has no periodicity. Therefore, we consider the equation in the complex plane directly. Denote by
$$g(\lambda)\triangleq-\dfrac{1}{2}[{\rm e}^{\tau\lambda}+{\rm e}^{(\tau-2)\lambda}].$$
For $a,b\in\mathbb{Z}\backslash\{0\}, a<b$, denote by the set
$$
\mathscr{C}_{a,b}\triangleq\left\{g(\lambda) \Bigg| \lambda\in\mathbb{C}_0, a\pi\leq{\rm Im}\lambda\leq b \pi \right\} \cap \mathbb{R}.
$$

\begin{figure}[!htb]\centering
 \includegraphics[width=0.48\textwidth]{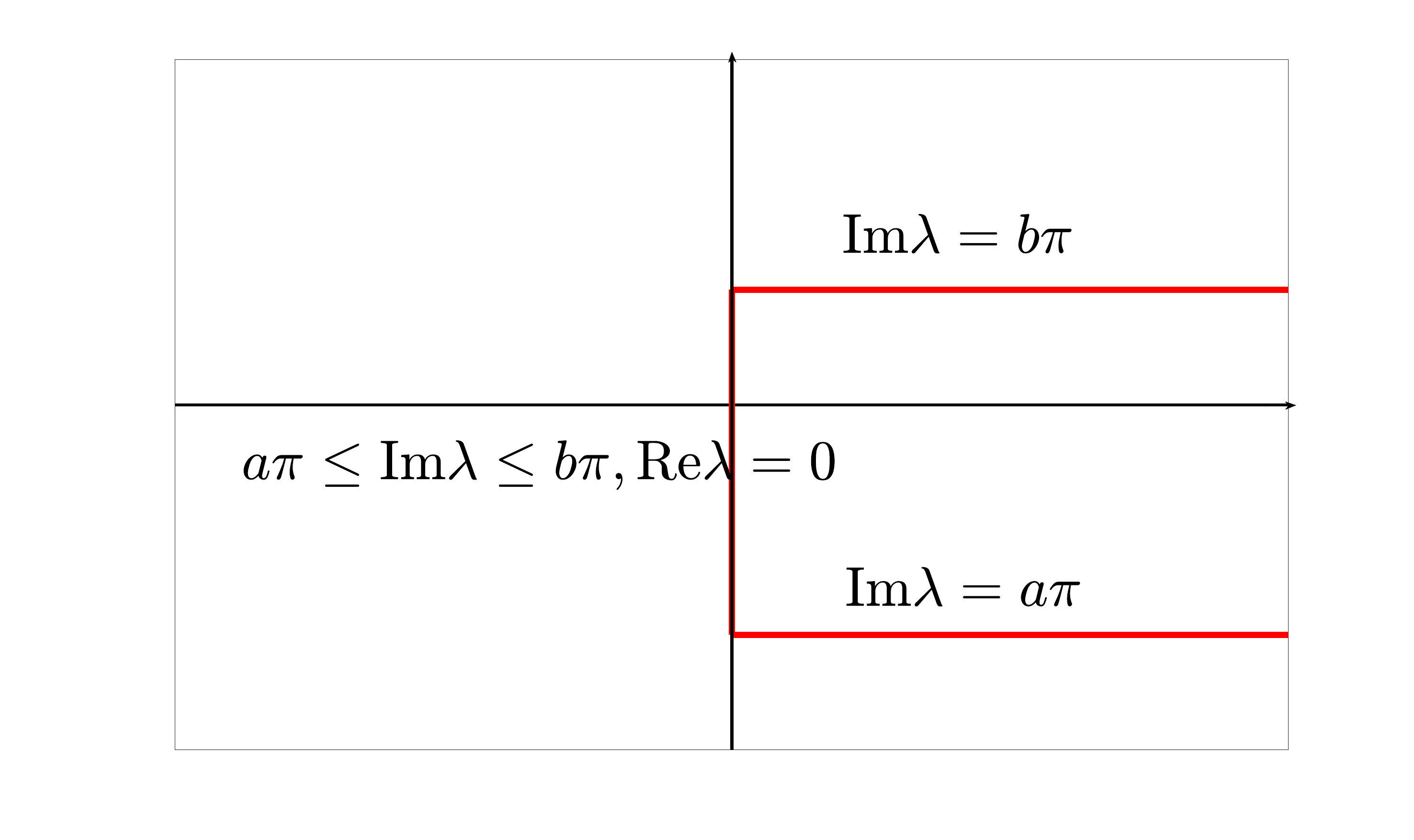}
\caption{The vertical strip $\left\{\lambda \Bigg| \lambda\in\mathbb{C}_+, a\pi\leq{\rm Im}\lambda\leq b \pi \right\}$.}\label{fig16}
\end{figure}

For $c\notin\mathscr{C}_{a,b}$, we denote by a function $M_{a,b}(c)$, which is the number of root, counted by multiplicity, of the equation $g(\lambda)=c$ located in the vertical strip $\left\{\lambda \Bigg| \lambda\in\mathbb{C}_+, a\pi\leq{\rm Im}\lambda\leq b \pi \right\}$ (See Fig. \ref{fig16}).

We only need to prove that there exist a pair of $a,b\in\mathbb{Z}\backslash\{0\}, a<b$ such that $M_{a,b}(c)$ is not zero for all $c\in\mathbb{R}$. The idea of the proof is quite similar to that in the last subsection. $M_{a,b}(c)$ is the number of root branches for different values of $c$. We will prove that $g(\lambda)\notin\mathbb{R}$ when ${\rm Im}\lambda=a\pi, b\pi$ for $a,b\in\mathbb{Z}\backslash{0}$. This indicates there is no root branch disappear or emerge at the boundary $\Big\{\lambda \Big|{\rm Im}\lambda=a\pi, b\pi\Big\}$. Each time when $c$ moves across some critical value $c\in\mathscr{C}_{a,b}$, a root branch disappear or emerge at the imaginary axis. We use implicit theorem to investigate whether a root branch disappears of emerges. When $\tau<1$, we move $c$ from $\pm\infty$ to $0$. When $c$ is at $\pm\infty$, we prove that there are at least one root branch by Lemma \ref{lemma8}. Each time when $c$ comes arcoss a critical value $c$, we prove that a new branch emerges  from the imaginary axis. Then we prove that there is no stability region. When $\tau>1$, we move $c$ from $0$ to $\pm\infty$. We prove that each time when $c$ comes arcoss a critical value $c\neq 0$, a new branch emeges from the imaginary axis. Thus, we only need to investigate the root branch near $c=0$. We prove that no matter $c$ moves from $0$ to $0+$ or $0-$, there is at least one root branch emerging and thus prove that there is no stability region for $c$.

Case A: $\tau<1$. According to Lemma \ref{lemma11},
we get that
$
{\rm Sgn}[{\rm Re}\lambda'(c_*)]={\rm Sgn}(\tau-1){\rm Sgn}(c_*)=-{\rm Sgn}(c_*).
$
Similar to the argument in last subsection, we obtain that $M_{a,b}(c)$ is monotonically increasing on $\mathbb{R}_{-}\backslash\mathscr{C}_{a,b}$ and monotonically decreasing on $\mathbb{R}_+\backslash\mathscr{C}_{a,b}$.  Lemma \ref{lemma8} shows that if we take $(b-a)\tau>2$, $M_{a,b}(c)$ is a nonzero constant when $|c|>1$. Therefore, $M_{a,b}(c)$ is not zero on $\mathbb{R}\backslash\mathscr{C}_{a,b}$.

Case B: $\tau>1$. We obtain that  $M_{a,b}(c)$ is monotonically increasing on $\mathbb{R}_{+}\backslash\mathscr{C}_{a,b}$ and monotonically decreasing on $\mathbb{R}_{-}\backslash\mathscr{C}_{a,b}$. We need to investigate $M_{a,b}(c)$ near $c=0$.  According to Lemma \ref{lemma11}, there exist $b-a-1$ different implicit functions $\lambda_k(c)$ such that $c=g(\lambda_k(c)), \lambda_k(0)=\lambda_k={\rm i}(k+\dfrac{1}{2})\pi$ , $a+1\leq k\leq b-1, k\in\mathbb{Z}$.  Furthermore,
$
{\rm Sgn}[{\rm Re}\lambda_k'(c_*)]=-{\rm Sgn}\cos[\tau(k+\dfrac{1}{2})\pi].
$
According to Lemma \ref{lemma12}, there exists $k,j\in\mathbb{N}^*$ such that
$\cos[\tau(j+\dfrac{1}{2})\pi]>0, \cos[\tau(l+\dfrac{1}{2})\pi]<0.$
Therefore, we take $a,b\in\mathbb{Z}\backslash\{0\}$ such that $a\leq\min\{l,j\}-1, b\geq \max\{l,j\}+1$. Then we get that
$
{\rm Sgn}[{\rm Re}\lambda_j'(c_*)]<0, {\rm Sgn}[{\rm Re}\lambda_l'(c_*)]>0.
$
Therefore, we get that $\lambda_l(c)\in\mathbb{C}_+$ for $c\in(0,\epsilon)$ and  $\lambda_j(c)\in\mathbb{C}_+$ for $c\in(-\epsilon,0)$. This leads to that $M_{a,b}(c)$ is not zero for $c\to 0+$ and $c\to 0-$. Thus, $M_{a,b}(c)$ is nonzero on $\mathbb{R}\backslash\mathscr{C}$.

From above all, the sufficient and necessary condition for the stability region for the parameter $c, \tau$ can be summarized as follows.
\begin{equation} \label{condition}
\left\{\begin{array}{l}
   c\in(-\sin [\dfrac{\pi}{2(\tau-1)}],0), \quad ~\mbox{for}~~\tau=4l-2, l\in \mathbb{N}^*,
      \\  c\in(0,\sin [\dfrac{\pi}{2(\tau-1)}]), \quad ~\mbox{for}~~\tau=4l, l\in\mathbb{N}^*.
\end{array}\right.
\end{equation}

It follows from Lemma \ref{Lm20199101844} that we get the following theorem.
\begin{theorem}\label{Th20189231016}
Suppose that $c $ and $\tau$  satisfy condition \eqref{condition}. Then, for any  initial state $(z( \cdot, 0),z_t(\cdot,0), w(\cdot, 0) )\in  \mathcal{X}$,
the closed-loop  system  \eqref{20189202158} with $c_1=c_2=c$
 admits a unique solution
$(z( \cdot, t),z_t(\cdot,t), w(\cdot, t) )\in C([0,\infty);\mathcal{X})$ which satisfies
 \begin{equation}\label{20189231017}
 \begin{array}{l}
 \|(z( \cdot, t),z_t(\cdot,t), w(\cdot, t) )\|_{ \mathcal{X} }\leq L_1{\rm e}^{-\omega_1  t}  \|(z( \cdot, 0),z_t(\cdot,0), w(\cdot, 0) )\|_{ \mathcal{X} } ,
\end{array}
\end{equation}
where    $L_1 $ and $\omega_1 $ are positive constants  independent of time and initial state.
\end{theorem}

\begin{remark}\label{remarkwangdelay2011}
We talk about the situation when $c_1=c_2=c$. If we take $c_1=0$, then system  \eqref{2018923949}  becomes \eqref{20199111503}, which has been fully investigated in \cite{Wangdelay2011}. In \cite{Wangdelay2011}, authors proved that when $\tau$ is an even number, there is not empty stability region for $c_2\in\mathbb{R}$. However, they did not prove the necessity of it and  could not provide the general formula for the stability region of $c_2$ for different values of $\tau$. They did not talk about the situation in which $\tau$ is irrational. If we use the same method employed in this section, we could prove that if and only if $\tau>0$ is an even number, there is  a stability region for $c_2$. Moreover, when $\tau$ is even number, the stability region can be summarized as follows.
\begin{equation} \label{condition2}
\left\{\begin{array}{l}
   c_2\in(-\tan (\frac{\pi}{2\tau}),0), \quad ~\mbox{for}~~\tau=4l-2, l\in \mathbb{N}^*,
      \\  c_2\in(0,\tan (\frac{\pi}{2\tau})), \quad ~\mbox{for}~~\tau=4l, l\in\mathbb{N}^*.
\end{array}\right.
\end{equation}
Obviously, the shrink of the stability region as $\tau$ is increasing can be explicitly obtained by $$\lim_{\tau\to+\infty}\tan(\dfrac{\pi}{2\tau})=0,$$
{which improves the results in \cite[Section 6]{Wangdelay2011}.  }
\end{remark}

\begin{remark}
For the situation $\tau>0$ is irrational, we can get the conclusion that \eqref{cha} has unstable roots from  \cite[Page 287, Page 288, Eq. (6.11)]{Hale}, which can be described as follows.

If $r_1,r_2>0$ are rationally independent, the sufficient and necessary condition for all the roots of the characteristic equation $1=a_1{\rm e}^{-\lambda r_1}+a_2{\rm e}^{-\lambda r_2}+a_3{\rm e}^{-\lambda (r_1+r_2)}$ lie in $\mathbb{C}_{-}$ is $1+a_1>|a_2+a_3|,1-a_1>|a_2-a_3|$.

Note that the characteristic equation \eqref{cha} can be written as $$1=-{\rm e}^{-2\lambda}-(c_1+c_2){\rm e}^{-\tau\lambda}-(c_1-c_2){\rm e}^{-(2+\tau)\lambda}.$$ 
By taking $a_1=-1, a_2=-(c_1+c_2), a_3=-(c_1-c_2)$ into  $1+a_1>|a_2+a_3|,1-a_1>|a_2-a_3|$, we easily get that the stable region for $c_1,c_2$ is empty.

\end{remark}

\section{Robustness to a small perturbation in time delay in low frequencies}\label{section4.5}

Enlightened by \cite{Datko86,Datko88,Datko93,Wangdelay2011}, we know that the feedback loop is not robust to a small perturbation in time delay. Authors in those literatures demonstrate the lack of robustness by giving exact expressions for eigenvalues for a special sequence of delay perturbations (See \cite[Theorem 7.2]{Wangdelay2011}, \cite[Lemma 2]{Datko86}) and thus in \cite[Page 5, Remark]{Datko86}, the author guessed that a small perturbation of $\epsilon$ in time delay will excite a high frequency mode~(i.e., a mode with frequency $\approx\mathcal{O}(\dfrac{1}{ |\epsilon| })$ as $\epsilon\to0$).  In this section, we will verify this judgement by spectral analysis for system \eqref{20189202158}. For the sake of simplicity, we only consider the situation $c_1=c_2=c$.

Firstly, we discuss about the robustness for $\tau=0$. When $\tau=0$, the stability region for $c$ is $(-\infty,-1)\cup(0,+\infty)$. For the sake of simplicity, we consider the robustness when $c>0$.

\begin{theorem}\label{robust4}
 Consider Eq.~\eqref{202302042155} when  $\tau_\epsilon=\epsilon~(\epsilon>0)$ and $c>0$. Denote by $\lambda_\epsilon\triangleq\inf\{|{\rm Im}\lambda|\Big|\lambda~\mbox{is a root of }~Eq.~\eqref{202302042155}~\mbox{located in}~\mathbb{C}_+\cup\mathbb{C}_0\}$, then $\lambda_\epsilon=\mathcal{O}(\dfrac{1}{\epsilon})$ as $\epsilon\to 0+$. This implies that there exists a positive constant $C_1$ independent of $\epsilon$, such that Eq. \eqref{202302042155} has no roots located in $\Big\{\lambda\in\mathbb{C}_0\cup\mathbb{C}_+\Big||{\rm Im}\lambda|<\dfrac{C_1}{\epsilon}\Big\}$ provided $\epsilon$ is sufficiently small.
 \end{theorem}
Proof We will prove that there exists two positive constants $C_1,C_2$~(independent of $\epsilon$) such that Eq. \eqref{202302042155} has no roots located in $\Big\{\lambda\in\mathbb{C}_0\cup\mathbb{C}_+\Big||{\rm Im}\lambda|<\dfrac{C_1}{\epsilon}\Big\}$ while Eq. \eqref{202302042155} has at least one root located in $\Big\{\lambda\in\mathbb{C}_0\cup\mathbb{C}_+\Big||{\rm Im}\lambda|<\dfrac{C_2}{\epsilon}\Big\}$. The proof will be divided into two parts. For the first part, we prove that Eq. \eqref{202302042155} has no roots located in $\Big\{\lambda\in\mathbb{C}_0\cup\mathbb{C}_+\Big||{\rm Im}\lambda|<\dfrac{C_1}{\epsilon}\Big\}$. For the second part, we prove that Eq. \eqref{202302042155} has at least one root located in $\Big\{\lambda\in\mathbb{C}_0\cup\mathbb{C}_+\Big||{\rm Im}\lambda|<\dfrac{C_2}{\epsilon}\Big\}$.

For the first part. We choose $C_1=\dfrac{\pi}{2}$.  Suppose that $\hat{\lambda}_\epsilon=p+{\rm i}q, p\geq 0, |q|<\dfrac{\pi}{2\epsilon}$ is a root of \eqref{202302042155} when $\tau=\epsilon, c>0$. We write \eqref{202302042155} as
$${\rm e}^{2\hat{\lambda}_\epsilon}(1+2c{\rm e}^{-\epsilon\hat{\lambda}_\epsilon})=-1,$$which, by taking absolute value of both sides, leads to
\begin{equation}\label{202302151106}
{\rm e}^{4p}(1+4c^2{\rm e}^{-2\epsilon p}+4c{\rm e}^{-\epsilon p}\cos {\epsilon q})=1.
\end{equation}
Since $|q|<\dfrac{\pi}{2\epsilon}$, $\cos{\epsilon q}>0$. Then we have
$1+4c^2{\rm e}^{-2\tau p}+4c{\rm e}^{-\tau p}\cos {\epsilon q}>1$, and ${\rm e}^{4b}\geq 1$. This contradicts \eqref{202302151106} and completes the proof of the first part.

For the second part, we choose $C_2>\pi$. Denote by $S_\epsilon$ as the smallest integer number such that $S_\epsilon>\dfrac{1}{\epsilon}$. Since $\tau=\epsilon<1$, by the proof in subsection \ref{subsection32}~(when $(b-a)\tau>2$, $M_{a,b}(c)$ is nonzero), we know that $M_{-S_\epsilon,S_\epsilon}(c)$~(defined in \ref{subsection32}) is nonzero when $\tau$ is irrational. This implies that Eq. $\eqref{202302042155}$ has at least one root located in $\Big\{\lambda\in\mathbb{C}_0\cup\mathbb{C}_+\Big||{\rm Im}\lambda|<S_\epsilon\pi\Big\}\subseteq \Big\{\lambda\in\mathbb{C}_0\cup\mathbb{C}_+\Big||{\rm Im}\lambda|<\dfrac{C_2}{\epsilon}\Big\}$. For rational $\tau=\epsilon>0$ efficiently small, we can use the similar idea of Lemma \ref{lemma7},\ref{lemma8} and \ref{lemma11} to prove that Eq. $\eqref{202302042155}$ has at least one root located in $\Big\{\lambda\in\mathbb{C}_0\cup\mathbb{C}_+\Big||{\rm Im}\lambda|<S_\epsilon\pi\Big\}$. We put the details of the proof in Appendix, Lemma \ref{lemma16}. This completes the proof.

Secondly, we discuss about the robustness for $\tau=2l, l\in\mathbb{N}^*$.

\begin{theorem}\label{robust3}
Consider Eq.~\eqref{202302042155} when $\tau_{\epsilon}=2l+\epsilon~(\epsilon\in\mathbb{R},l\in\mathbb{N}^*)$ and $c$ for which $2l,c$ satisfy Condition \eqref{condition}. Denote by $\lambda_\epsilon\triangleq\inf\{|{\rm Im}\lambda|\Big|\lambda~\mbox{is a root of }~Eq.~\eqref{202302042155}~\mbox{located in}~\mathbb{C}_+\cup\mathbb{C}_0\}$, then $\lambda_\epsilon=\mathcal{O}(\dfrac{1}{|\epsilon|})$ as $|\epsilon|\to 0$. This implies that there exsits a positive constant $C_1$ independent of $\epsilon$, such that Eq. \eqref{202302042155} has no roots located in $\Big\{\lambda\in\mathbb{C}_0\cup\mathbb{C}_+\Big||{\rm Im}\lambda|<\dfrac{C_1}{|\epsilon|}\Big\}$ provided $|\epsilon|$ is sufficiently small.

\end{theorem}

Proof
We will prove that there exists two positive constants $C_1,C_2$~(independent of $\epsilon$) such that Eq. \eqref{202302042155} has no roots located in $\Big\{\lambda\in\mathbb{C}_0\cup\mathbb{C}_+\Big||{\rm Im}\lambda|<\dfrac{C_1}{|\epsilon|}\Big\}$ while Eq. \eqref{202302042155} has at least one root located in $\Big\{\lambda\in\mathbb{C}_0\cup\mathbb{C}_+\Big||{\rm Im}\lambda|<\dfrac{C_2}{|\epsilon|}\Big\}$. The proof will be divide into two parts. For the first part, we prove that Eq. \eqref{202302042155} has no roots located in $\Big\{\lambda\in\mathbb{C}_0\cup\mathbb{C}_+\Big||{\rm Im}\lambda|<\dfrac{C_1}{|\epsilon|}\Big\}$. For the second part, we prove that Eq. \eqref{202302042155} has at least one root located in $\Big\{\lambda\in\mathbb{C}_0\cup\mathbb{C}_+\Big||{\rm Im}\lambda|<\dfrac{C_2}{|\epsilon|}\Big\}$.

Taking $\tau_{\epsilon}=2l+\epsilon$ and $c$ in Eq. \eqref{202302042155} leads to the equation $h_\epsilon(\lambda)=c$, where
$$
h_\epsilon(\lambda)\triangleq -\dfrac{1}{2}{\rm e}^{\epsilon\lambda}{\rm e}^{2l\lambda}(1+{\rm e}^{-2\lambda}).
$$
Since $0<|c|<\sin[\dfrac{\pi}{2(2l-1)}]$, there exists $\tilde{c}\in(0,1)$ such that $|c|=\sin[\dfrac{\tilde{c}\pi}{2(2l-1)}]$. We choose $C_1\triangleq\dfrac{(1-\tilde{c})\pi}{2}>0, C_2\triangleq\dfrac{\pi}{2}$ and $s_\epsilon$ as the smallest integer number such that $s_\epsilon\pi>\dfrac{C_1}{|\epsilon|}$, while $S_\epsilon$ as the largest integer number such that  $S_\epsilon\pi<\dfrac{C_2}{|\epsilon|}$.

 Firstly, we prove that  for sufficiently small $|\epsilon|$, the equation $h_{\epsilon}(\lambda)=c$ has no roots lie in the vertical stripe $\left\{\lambda \Bigg| {\rm Re\lambda}\geq 0,|{\rm Im}\lambda|\leq s_\epsilon\pi \right\}$.
 For a fixed $\epsilon\neq 0$ and $\delta\in[-|\epsilon|,|\epsilon|]$, denote by $Q_\epsilon(\delta)$ as the number of roots  located in the vertial stripe $\left\{\lambda \Bigg| {\rm Re\lambda}\geq 0,|{\rm Im}\lambda|\leq s_\epsilon\pi \right\}$ for the equation $h_\delta(\lambda)=c$. By argument principle, we obtain that
$$
Q_\epsilon(\delta)=\dfrac{1}{2\pi{\rm i}}\int_{P_{\epsilon,R}}\dfrac{h'_{\delta}(\lambda)}{h_{\delta}(\lambda)-c}{\rm d}\lambda.
$$

\begin{figure}[!htb]
  \centering
    \subfigure[]{\includegraphics[width=0.49\textwidth]{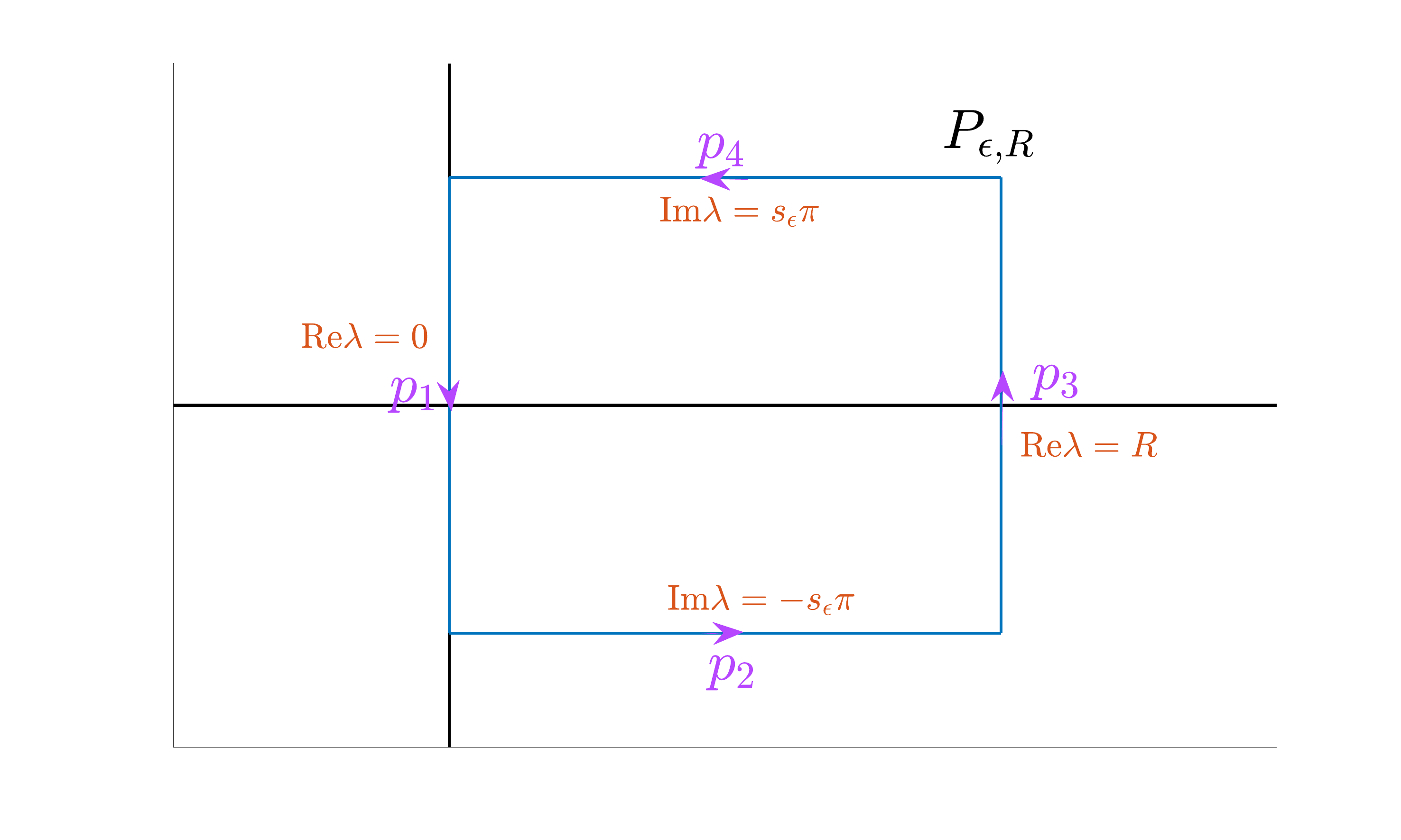}\label{fig17a}}
    \subfigure[]{\includegraphics[width=0.49\textwidth]{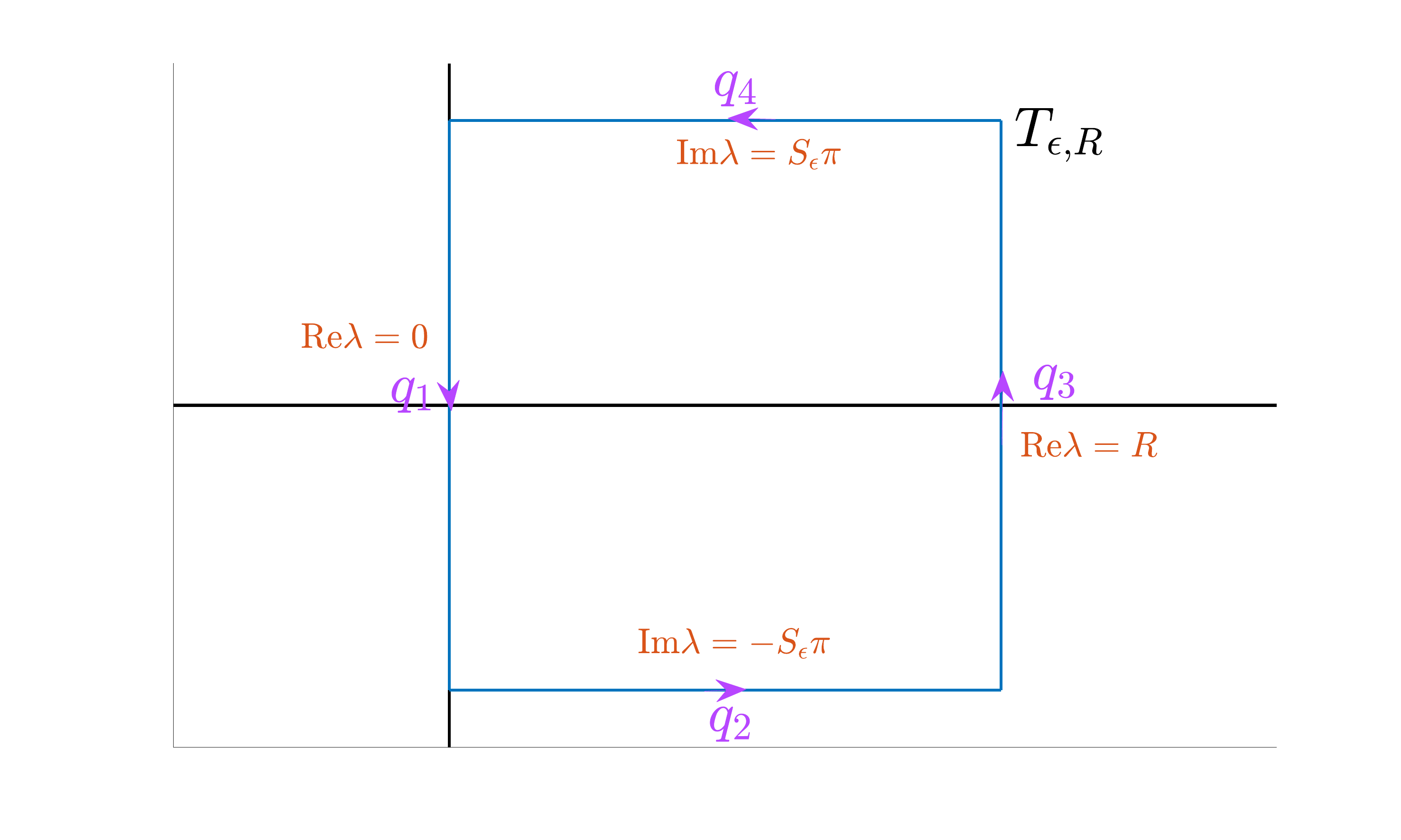}\label{fig17b}}
   \caption{The rectangle contour $P_{\epsilon,R}$  and $T_{\epsilon,R}$. }
    \label{fig17}
   \end{figure}

Here, $P_{\epsilon,R}$ denotes the rectangle contour $p_1\cup p_2\cup p_3\cup p_4$, where
$$
p_1\triangleq\{\lambda={\rm i}\beta, \beta:s_\epsilon\pi\to -s_\epsilon\pi\}, p_2\triangleq\{\lambda=u-{\rm i}s_\epsilon\pi, u:0\to R \},
$$
$$
p_3\triangleq\{\lambda=R+{\rm i}\beta, \beta:-s_\epsilon\pi\to s_\epsilon\pi\}, p_4\triangleq\{\lambda=u+{\rm i}s_\epsilon\pi, u:R\to 0 \},
$$
and $R>0$ is a sufficiently large number such that $|h_\delta(\lambda)|>1>|c|$ for ${\rm Re}\lambda\geq R$ and $\delta\in[-|\epsilon|,|\epsilon|]$~(See Fig. \ref{fig17a}). The existence of $R>0$ is guaranteed by  $\lim_{{\rm Re}\lambda\to+\infty}|h_\delta(\lambda)|=+\infty$, where the limit is taken uniformly with respect to $\delta\in[-|\epsilon|,|\epsilon|]$.

Condition \eqref{condition} ensures that $h_0(\lambda)=c$ has no roots in $\mathbb{C}_0\cup\mathbb{C}_+$, which indicates that $Q_\epsilon(0)=0$. Similar to the proof in Section \ref{section3}, we hope to prove that $Q_\epsilon(\delta)$ is continuous with respect to $\delta\in[-|\epsilon|,|\epsilon|]$. For this, we just need to prove that $h_\delta(\lambda)\neq c$ on the rectangle $P_{\epsilon,R}$. For $\lambda\in p_3$, $|h_\delta(\lambda)|>1>|c|$, which indicates $h_\delta(\lambda)\neq c$ on the line $h_3$. For $\lambda\in p_2\cup p_4$, $|{\rm Im}\lambda|=s_\epsilon\pi$, which indicates that ${\rm e}^{2l\lambda}(1+{\rm e}^{-2\lambda})\in\mathbb{R}\backslash\{0\}$. Furthermore, when $\delta\neq 0$, the argument of ${\rm e}^{\delta\lambda}$  for $\lambda\in p_2\cup p_4$ can be estimated as
$$
0<|{\rm arg}~({\rm e}^{\delta\lambda})|=|\delta{\rm Im}\lambda|\leq |\epsilon|s_\epsilon\pi\leq C_1+|\epsilon|\pi<(\dfrac{1-\tilde{c}}{2}+|\epsilon|)\pi<\pi.
$$
Therefore, we obtain ${\rm e}^{\delta\lambda}\notin\mathbb{R}$. This implies that $h_\delta(\lambda)\notin\mathbb{R}$ and thus $h_\delta(\lambda)\neq c$ for $\lambda\in p_2\cup p_4$. For $\delta=0$, we use Condition \eqref{condition} to ensure that $h_{0}(\lambda)\neq c$ on the line $p_2\cup p_4$. Therefore, we have proved that $h_\delta(\lambda)\neq c$ for $\lambda\in p_2\cup p_4,\delta\in[-|\epsilon|,|\epsilon|]$.

The most complicated part is to prove $h_{\delta}(\lambda)\neq c$ for $\lambda\in p_1$. We put this in the Appendix, Lemma \ref{lemma13} and completes the first part of proof.

Secondly, we prove that  for sufficiently small $|\epsilon|$, the equation $h_{\epsilon}(\lambda)=c$ has at least one root lie in the vertical stripe $\left\{\lambda \Bigg| {\rm Re\lambda}\geq 0,|{\rm Im}\lambda|\leq S_\epsilon\pi \right\}$. We only prove the case $\epsilon>0$. The case $\epsilon<0$ is similar so we omit it.

For a fixed $\epsilon>0$ and $\delta\in[-\epsilon,\epsilon]$, denote by $H_\epsilon(\delta)$ as the number of roots  located in the vertial stripe $\left\{\lambda \Bigg| {\rm Re\lambda}\geq 0,|{\rm Im}\lambda|\leq S_\epsilon\pi \right\}$ for the equation $h_\delta(\lambda)=c$. By argument principle, we obtain that
$$
H_\epsilon(\delta)=\dfrac{1}{2\pi{\rm i}}\int_{T_{\epsilon,R}}\dfrac{h'_{\delta}(\lambda)}{h_{\delta}(\lambda)-c}{\rm d}\lambda.
$$
Here,  $T_{\epsilon,R}$ denotes the rectangle contour $q_1\cup q_2\cup q_3\cup q_4$, where
$$
q_1\triangleq\{\lambda={\rm i}\beta, \beta:S_\epsilon\pi\to -S_\epsilon\pi\}, q_2\triangleq\{\lambda=u-{\rm i}S_\epsilon\pi, u:0\to R \},
$$
$$
q_3\triangleq\{\lambda=R+{\rm i}\beta, \beta:-S_\epsilon\pi\to S_\epsilon\pi\}, q_4\triangleq\{\lambda=u+{\rm i}S_\epsilon\pi, u:R\to 0 \},
$$
while $R>0$ is sufficiently large~(See Fig. \ref{fig17b}). Similarly to the first part, we get that $h_\delta(\lambda)\neq c$ on the line $q_2\cup q_3\cup q_4$. We only need to investigate $h_\delta(\lambda)$ on $q_1$. Akin to the proof in Section \ref{section3}, denote by a set
$$
\mathscr{F}_{\epsilon}\triangleq\Big\{ \delta\in[0,\epsilon]\Big |   \mbox{there exists}~~\beta\in\mathbb{R}~~\mbox{such that}~~h_\delta({\rm i}\beta)=c, |\beta|\leq S_\epsilon\pi   \Big\}.
$$
Thus, $\mathscr{F}_{\epsilon}$ is the set of the critical value for $H_\epsilon(\delta)$.  We also need to investigate whether $H_\epsilon(\delta+)-H_\epsilon(\delta-)$ is positive or negative by the Implicit Theorem for each $\delta\in\mathscr{F}_{\epsilon}$. We prove that $\mathscr{F}_{\epsilon}$ is nonempty~(Lemma \ref{lemma14}) and furthermore for each $\delta\in\mathscr{F}_{\epsilon}$,  $H_\epsilon(\delta+)-H_\epsilon(\delta-)>0$~(Lemma \ref{lemma15}). We put all the details in the Appendix, Lemma \ref{lemma14} and Lemma \ref{lemma15} and complete the proof of the second part.


%

\begin{remark}\label{robustremark}
As noted by \cite{Wangdelay2011, Datko86, Logemann}, {\color{black} the feedback stabilizer of wave equation    usually shows lack of robustness to a small  delay perturbations. 
}Despite this non-robustness, numerical experiments often demonstrate an absence of the destabilizing effect when a small perturbation is added to the time delay. 
It is important to note that numerical experiments often neglect high frequency modes. Theorem \ref{robust4} and \ref{robust3} confirm that when a small perturbation is added to the time delay, no roots are found in low frequencies.\end{remark}

\section{Appendix}

\begin{lemma}\label{lemma1}
The set $\mathscr{E}$ can be characterized as

$$
\mathscr{E}_{m,n}=\Big\{(-1)^{l}\cot{\dfrac{n(2l+1)\pi}{2m}} \Big|   
l=0,1,2,\cdots, 2m-1          \Big\}\cup \Big\{0 \Big\}.
$$

\end{lemma}
Proof For an element $k\in\mathscr{E}$, there exists a complex number $z={\rm e}^{{\rm i}\theta}, \theta\in [0,2\pi)$, such that
$c=f(z)$. Thus,

By taking $z={\rm e}^{{\rm i}\theta}, \theta\in[0,2\pi)$ into $f(z,k)=0$, we obtain that
Case A: $ k=0,   z^{2n}+1=0$

Case B: $\cos(n\theta)=-k{\rm e}^{{\rm i}m\theta}{\rm i}\sin(n\theta)$

$$
c=-\dfrac{1}{2}\bigg\{\cos[(2n-m)\theta]+\cos (m\theta)\bigg\}
-\dfrac{\rm i}{2}\bigg\{\sin[(2n-m)\theta]-\sin (m\theta)\bigg\}.
$$
Since $\sin[(2n-m)\theta]-\sin (m\theta)=0$, we get that $\theta=\dfrac{k\pi}{|m-n|}$, $k=0,1,2,\cdots, 2|m-n|-1$ and in this case, $s=-\cos(m\theta)$
or $\theta=\dfrac{(2k+1)\pi}{2n}, k=0,1,2,\cdots, 2n-1$ and in this case, $c=0$. This completes the proof.

\begin{lemma}\label{lemma2}
For each $c_*=f(z_*)\in\mathscr{E}\backslash\{0\}, |z_*|=1$, there exists an implicit function $z(c)$ such that $z(c_*)=z_*$ and $c=f(z(c))$ for each $c\in(c_*-\epsilon, c_*+\epsilon)$. Here, $\epsilon>0$ is sufficiently small. Suppose that $z(c)=r(c){\rm e}^{{\rm i}\theta(c)}$, where $r(c)$ and $\theta(c)$ represent the absolute value and argument value of the function $z(c)$, respectively. Then we have
$$
{\rm Sgn}[r'(c_*)]={\rm Sgn}(n-m){\rm Sgn}(c_*).
$$
For $c_*=0$, there exist $2n$ different implicit functions $z_k(c)$ such that $z_k(0)=z_k={\rm e}^{\rm i \theta_k}$ with $\theta_k=\dfrac{(2k+1)\pi}{2n}$, $k=0,1,2,\cdots, 2n-1$.  And we have
$$
{\rm Sgn}[r_k'(0)]={\rm Sgn}\cos[\dfrac{m(2k+1)\pi}{2n}].
$$
\end{lemma}
Proof Since $0=f(z(k),k)$, Taking derivatives with respect to $k$ of both sides of   $f(z(k),k)=0$   leads to
$z'(k)=-\dfrac{f_k}{f_z}$
Taking derivatives with respect to variable $k$ of both sides  of   $z=r{\rm e}^{\rm i\theta}$    leads to that
$\dfrac{z'}{z}=\dfrac{r'}{r}+{\rm i}\theta'$,
which further indicates that ${\rm Sgn}[r'(c)]={\rm Sgn}\{{\rm Re}[\dfrac{1}{f'(z)z}]\}$
For $c_*=f(z_*)\in\mathscr{E}\backslash\{0\}$, then $z_*=-\cos(m\theta)$ with $\theta=\dfrac{k\pi}{|m-n|}$, $k=0,1,2,\cdots, 2|m-n|-1$. We compute that ${\rm Re} [f'(z_*)z_*]=-(n-m)\cos(m\theta)=(n-m)c_*\neq 0$, which ensures the existence of the implicit function $z(c)$ near each nonzero $c_*\in\mathscr{E}$.
Furthermore, we get ${\rm Sgn}[r'(c_*)]={\rm Sgn}(n-m){\rm Sgn}(c_*)$.

 For $k_*=0$, the equation $f(z,k)=0$ has $2n$ different roots $z_l={\rm e}^{\rm i \theta_l}$ with $\theta_l=\dfrac{(2l+1)\pi}{2n}$, $k=0,1,2,\cdots, 2n-1$.

 Since $z_l^{2n}=-1$, we compute that$f'(z_k)z_k=nz_k^{-m}\neq 0$, which ensures the existence of each implicit function $z_k(c)$ and further implies that
${\rm Sgn}[r_k'(0)]={\rm Sgn}\cos[\dfrac{m(2k+1)\pi}{2n}]$. This completes the proof.

\begin{lemma}\label{lemma3}
Denote by the set
$$
\mathscr{B}=\left\{  \cos[\dfrac{m(2k+1)\pi}{2n}]  \Bigg|k=0,1,2,\cdots, 2n-1 \right\}.
$$
If $m,n$ are coprime positive integers and $m>n\geq 2$, then $\mathscr{B}\cap \mathbb{R}_+\neq \emptyset$ and $\mathscr{B}\cap \mathbb{R}_-\neq \emptyset$.
\end{lemma}
Proof Case I: $m$ is an odd number. Since $m$ and $2n$ are coprime, for any integer number $x$, there exists integer number $p\in\{0,1,2,\cdots, 2n-1\}$ and $q$ such that $mp-2nq=x$. Thus, we can find integer number $p_1,q_1,p_2,q_2$ with $p_1,p_2\in\{0,1,2,\cdots, 2n-1\}$ such that
$mp_1-2nq_1=\dfrac{1-m}{2}$, $mp_2-2nq_2=n-\dfrac{1+m}{2}$,
which implies that
 $\cos[\dfrac{m(2p_1+1)\pi}{2n}]=\cos(\dfrac{\pi}{2n})>0$ and $\cos[\dfrac{m(2p_2+1)\pi}{2n}]=-\cos(\dfrac{\pi}{2n})<0$.

Case II: $m$ is an even number and $n\geq 3$ is an odd number. In this case, $\cos(\dfrac{s\pi}{2n})>0$ and $\cos[\dfrac{(2n-s)\pi}{2n}]<0$ for $s=0,2.$
We find $s_1,s_2\in\{0,2\}$ such that $(s_1-m)/4$ and $(2n-m-s_2)/4$ are both integers. Since $\dfrac{m}{2}$ and $n$ are coprime, there exists integer number $p_1,q_1,p_2,q_2$ with $p_1,p_2\in\{0,1,2,\cdots, n-1\}$ such that
$\dfrac{m}{2}p_1-nq_1=\dfrac{s_1-m}{4}$,
$
\dfrac{m}{2}p_2-nq_2=\dfrac{2n-s_2-m}{4},
$
which further implies that
$
\cos[\dfrac{m(2p_1+1)}{2n}\pi]=\cos(\dfrac{s_1\pi}{2n})>0,
$
$
\cos[\dfrac{m(2p_2+1)}{2n}\pi]=\cos(\dfrac{2n\pi-s_2\pi}{2n})<0.
$

\begin{lemma}\label{lemma4}
If $n=1,m=4s-2, s\in\mathbb{N}^*$ , the largest element in the set $\mathscr{E}\cap \mathbb{R}_{-}$ is $-\sin[\dfrac{\pi}{2(m-1)}]$. If $n=1, m=4s, s\in\mathbb{N}^*$, the smallest element in the set $\mathscr{E}\cap \mathbb{R}_+$ is $\sin[\dfrac{\pi}{2(m-1)}]$.
\end{lemma}
Proof We only talk about the situation $n=1, m=4s-2$, the another is similar. Since $m$ and $m-1$ are coprime, there exists integer number $p,q$ with $p\in\{0,1,2,\cdots,m-2\}$ such that
$
\dfrac{m}{2}p-(m-1)q=\dfrac{m-2}{4},
$
which leads to
$
\dfrac{pm\pi}{m-1}-2q\pi=\dfrac{\pi}{2}-\dfrac{\pi}{2(m-1)}.
$
Therefore, we get that
$
-\cos [\dfrac{pm\pi}{m-1}]=-\sin [\dfrac{\pi}{2(m-1)}].$
On the other hand, for any integer number $k$ and $l$, if $\dfrac{km\pi}{m-1}\neq \dfrac{2l+1}{2}\pi$, then
$$
|\dfrac{km\pi}{m-1}-\dfrac{2l+1}{2}\pi|=|\dfrac{2km-(l+1)(m-1)}{2(m-1)}\pi|\geq \dfrac{\pi}{2(m-1)},$$ which implies that the absolute value of any nonzero element in $\mathscr{E}$ is larger than $\sin [\dfrac{\pi}{2(m-1)}].$ This completes the proof.

\begin{lemma}\label{lemma5}
If $n=1$ and $m\geq 3$ is an odd number, then $r_k'(0)=0, r_k''(0)<0$ for $k=0,1$, where the function $r_k$ is defined in Lemma \ref{lemma2}.
\end{lemma}
Proof When $n=1$ and $m\geq 3$ is an odd number, we get that $\cos[\dfrac{m(2k+1)\pi}{2n}]=0$ for $k=0,1$. According to Lemma \ref{lemma2}, we get that  $r_k'(0)=0$ for $k=0,1$.

Now we turn to compute $r_k''(0)$.
Taking derivatives of both sides of $1=f'(z)z'(c)$ leads to that
$0=f''(z)z'(c)^2+f'(z)z''(c)$,
and we obtain that $z''(c)=-\dfrac{f''(z)z'(c)^2}{f'(z)}=-\dfrac{f''(z)}{f'(z)^3}$.
Taking second derivatives of both sides of $z=r{\rm e}^{\rm i\theta}$ leads to that
$
z''=r''{\rm e}^{\rm i\theta}-r{\rm e}^{\rm i\theta}\theta'^2+{\rm i}[2r'\theta'{\rm e}^{\rm i\theta}+r{\rm e}^{\rm i\theta}\theta''],
$
which divided by $z$ leads to that
$
-\dfrac{f''(z)}{zf'(z)^3}=\dfrac{z''}{z}=\dfrac{r''}{r}-(\theta')^2+{\rm i}[\dfrac{2r'\theta'+\theta''}{r}].
$
We take the implicit function $z_k$ at point $c=0$. In this case $r_k(0)=1$,we compute that
$r_k''(0)=\theta_k'(0)^2+{\rm Re}[-\dfrac{f''(z)}{zf'(z)^3}]|_{z=z_k}= \theta_k'(0)^2-(2m-1)$.
From $\dfrac{z'}{z}=\dfrac{r'}{r}+{\rm i}\theta'$, we obtain that
$
\theta'={\rm Im}[\dfrac{z'}{z}]={\rm Im}[\dfrac{1}{zf'(z)}]={\rm Im} z^m.
$
The last equality is taken at $c=0, z=z_k$. This indicates that $|\theta_k'(0)|\leq 1$ and we obtain that $r_k''(0)<0$. This completes the proof.

\begin{lemma}\label{lemma6}
For an interval $[k_1,k_2]$, if $[k_1,k_2]\cap \mathscr{E}_{m,n}=\emptyset$, $N(k)$ is a constant. 
\end{lemma}
Proof
By argument principle , for $k\notin\mathscr{E}_{m,n}$, we have
\begin{equation}\label{argument}
N(k)-n=\dfrac{1}{2\pi{\rm i}}\int_{|z|=1}\dfrac{ f_z(z,k)}{f(z,k)}{\rm d}z.
\end{equation}
For an interval $[k_1,k_2]$, if $[k_1,k_2]\cap \mathscr{E}_{m,n}=\emptyset$, the right side is continuous with respect to variable $k\in[k_1,k_2]$. Since $N(k)$ is integer, we obtain that $N(k)$ is constant on the interval $[k_1,k_2]$.

Furthermore, we have
$
\lim_{|c|\to+\infty}\dfrac{1}{2\pi{\rm i}}\int_{|z|=1}\dfrac{f'(z)}{f(z)-c}{\rm d}z=0
$
and $\mathscr{E}\cap \big\{c\big||c|>1\big\}=\emptyset$, thus we get that
$N(c)\equiv m$ when $|c|>1$.

\begin{lemma}\label{lemma7}
For an interval $[c_1,c_2]$, if $[c_1,c_2]\cap \mathscr{C}_{a,b}=\emptyset$, $M_{a,b}(c)$ is a constant.
\end{lemma}
Proof
Since $\lim_{{\rm Re}\lambda\to+\infty}|g(\lambda)|=+\infty$, for any sufficiently large $R>0$, we find a sufficiently large $N>0$ such that $|g(\lambda)|>R$ for ${\rm Re} \lambda>N$. By argument principle \cite{complexanalysis}, we obtain that for  $c\in (-R,R)\backslash\mathscr{C}_{a,b}$,
\begin{equation}\label{argument2}
M_{a,b}(c)=\dfrac{1}{2\pi{\rm i}}\int_{P_{a,b,N}}\dfrac{g'(\lambda)}{g(\lambda)-c}{\rm d}\lambda.
\end{equation}
Here, $P_{a,b,N}$ denotes the rectangle contour $h_1\cup h_2\cup h_3\cup h_4$, where
$$
h_1\triangleq\{\lambda={\rm i}\beta, \beta:b\pi\to a\pi \}, h_2\triangleq\{\lambda=s+a{\rm i}\pi, s:0\to N \},
$$
$$
h_3\triangleq\{\lambda=N+{\rm i}\beta, \beta:a\pi\to b \pi\}, h_4\triangleq\{\lambda=s+b{\rm i}\pi, s:N\to 0 \}.
$$

Since $a,b\in\mathbb{Z}\backslash\{0\}$, we get $1+{\rm e}^{-2\lambda}\in\mathbb{R}_{+}$ on the line $h_2, h_4$, which further indicates that $g(\lambda)=-\dfrac{1}{2}(1+{\rm e}^{-2\lambda}){\rm e}^{\lambda\tau}\notin\mathbb{R}$ on the line $h_2, h_4$. Note that $g(\lambda)\neq c$ on $h_1, h_3$. Therefore, for an interval $[c_1,c_2]\subseteq (-R, R)$, if $[c_1,c_2]\cap \mathscr{C}_{a,b}=\emptyset$, the right side of Eq. \eqref{argument2} is continuous with respect to variable $c\in[c_1,c_2]$. Since $M_{a,b}(c)$ is integer number, $M_{a,b}(c)$ is a constant on interval $[c_1,c_2]$.

\begin{lemma}\label{lemma8}
If $(b-a)\tau>2$,$\tau$ is irrational, $M_{a,b}(c)$ is a nonzero constant  on the interval $(-\infty,-1)$ and $(1,+\infty)$, respectively.
\end{lemma}
Proof
Obviously, $|c|\leq 1$ for $c\in\mathscr{C}_{a,b}$. By Lemma \ref{lemma7}, we obtain that $M_{a,b}(c)$ is a constant on the interval $(1,+\infty)$ and $(-\infty,-1)$, respectively. We only need to prove that there exist  $\lambda_{1,2}\in\mathscr{C}_{a,b}$ such that $g(\lambda_1)\in(1,+\infty)$ and $g(\lambda_2)\in(-\infty,-1)$.

Suppose that $\lambda=p+{\rm i}q, p>0,a\pi<q<b\pi$. Thus, we get that
$$
g(\lambda)=-\dfrac{1}{2}[{\rm e}^{\tau p}\cos(\tau q)+{\rm e}^{\tau p-2p}\cos(\tau q-2q)]-\dfrac{{\rm i}}{2}[{\rm e}^{\tau p}\sin(\tau q)+{\rm e}^{\tau p-2p}\sin(\tau q-2q)].
$$
When ${\rm e}^{\tau p}\sin(\tau q)+{\rm e}^{\tau p-2p}\sin(\tau q-2q)=0,$ we compute that
$
{\rm e}^{2p}=-\dfrac{\sin(\tau q-2q)}{\sin(\tau q)},g(\lambda)=\dfrac{{\rm e}^{\tau p}\sin(2q)}{2\sin(\tau q-2q)}.
$

Since $b\tau-a\tau>2$,  there exists  odd number $k_1$ and even number $k_2$ such that $\tau a< k_i<\tau b$ for $i=1,2$. 

 Consider $q^*_i=\dfrac{k_i\pi}{\tau}, i=1,2$, since $\tau$ is irrational, we obtain that $\sin[(\tau-2)q^*_i]\neq 0$. Thus, we can take two sequences of $q_{i}^n,p_{i}^n, i=1,2, n=1,2,\cdots$, such that $\lim_{n\to+\infty}q_{i}^n=q^*_i$ and $\lim_{n\to+\infty}p_{i}^n=+\infty$  while  ${\rm e}^{2p_{i}^n}=-\dfrac{\sin(\tau q_{i}^n-2q_{i}^n)}{\sin(\tau q_{i}^n)} $   for $i=1,2$.

 For these two sequences, we obtain that
 $\dfrac{\sin(2q_i^n)}{\sin(\tau q_i^n-2q_i^n)}=-\lim_{n\to+\infty}\dfrac{\sin(\tau q_i^n-2q_i^n-\tau q_i^n)}{\sin(\tau q_i^n-2q_i^n)}=(-1)^{k_i+1}.$
If we take $\lambda_i^n=p_i^n+{\rm i}q_i^n$, we obtain that $g(\lambda_1^n)>0$, $\lim_{n\to+\infty}g(\lambda_1^n)=+\infty$ and  $g(\lambda_2^n)<0$, $\lim_{n\to+\infty}g(\lambda_2^n)=-\infty$. This completes the proof.

\begin{lemma}\label{lemma11}
For each nonzero $c_*=g(\lambda_*)\in\mathscr{C}_{a,b}, \lambda_*\in\mathbb{C}_0, a\pi\leq{\rm Im}\lambda_*\leq b\pi$, there exists an implicit function $\lambda(c)$ such that $\lambda(c_*)=\lambda_*$ and $c=g(\lambda(c))$ for each $c\in(c_*-\epsilon, c_*+\epsilon)$. Here, $\epsilon>0$ is sufficiently small. Furthermore, we have
$$
{\rm Sgn}[{\rm Re}\lambda'(c_*)]={\rm Sgn}(\tau-1){\rm Sgn}(c_*).
$$
For $c_*=0$, there exist $b-a-1$ different implicit functions $\lambda_k(c)$ such that $c=g(\lambda_k(c)), \lambda_k(0)=\lambda_k={\rm i}(k+\dfrac{1}{2})\pi$ , $a+1\leq k\leq b-1, k\in\mathbb{Z}$.  And we have
$$
{\rm Sgn}[{\rm Re}\lambda'(c_*)]=-{\rm Sgn}\cos[\tau(k+\dfrac{1}{2})\pi].
$$
\end{lemma}
ProofFor $c_*=g(\lambda_*)\in\mathscr{C}_{a,b}\backslash\{0\}$, since
$
c_*=g(\lambda_*)=-\dfrac{1}{2}[{\rm e}^{\tau\lambda_*}+{\rm e}^{(\tau-2)\lambda_*}]$,
we get that
$
{\rm Im} ({\rm e}^{\tau\lambda_*})=-{\rm Im}({\rm e}^{(\tau-2)\lambda_*}).$
By
$
|{\rm e}^{\tau\lambda_*}|=|{\rm e}^{(\tau-2)\lambda_*}|=1,
$
we obtain that $ |{\rm Re} ({\rm e}^{\tau\lambda_*})|=|{\rm Re}({\rm e}^{(\tau-2)\lambda_*})|$. Since $c_*\neq 0$, we have $ {\rm Re} ({\rm e}^{\tau\lambda_*})={\rm Re}({\rm e}^{(\tau-2)\lambda_*})=-\dfrac{c_*}{2}.$
Thus we have~
${\rm Sgn}[{\rm Re}\lambda'(c_*)]={\rm Sgn}[{\rm Re}\dfrac{1}{g'(\lambda_*)}]={\rm Sgn}[{\rm Re}g'(\lambda_*)]
=-{\rm Sgn}\{{\rm Re}[\tau {\rm e}^{\tau\lambda^*}+(\tau-2){\rm e}^{(\tau-2)\lambda^*}]\}={\rm Sgn}(\tau-1){\rm Sgn}(c_*).
$

For $c_*=0$, $g(\lambda)=0$ has $b-a-1$ different roots $\lambda_k={\rm i}(k+\dfrac{1}{2})\pi$ , $a+1\leq k\leq b-1, k\in\mathbb{Z}$. By ${\rm e}^{2\lambda_k}=-1$,  we compute that
$
g'(\lambda_k)=-\dfrac{1}{2}[\tau{\rm e}^{\tau\lambda_k}+(\tau-2){\rm e}^{(\tau-2)\lambda_k}]
=-\dfrac{1}{2}[\tau{\rm e}^{\tau\lambda_k}-(\tau-2){\rm e}^{\tau\lambda_k}]=-{\rm e}^{\tau\lambda_k}.
$
Then we have
$
{\rm Sgn}[{\rm Re}\lambda'(c_*)]={\rm Sgn}[{\rm Re}\dfrac{1}{g'(\lambda_*)}]={\rm Sgn}[{\rm Re}g'(\lambda_*)]=-{\rm Sgn}\cos[\tau(k+\dfrac{1}{2})\pi].
$
This completes the proof.

\begin{lemma}\label{lemma12}
There exist  $j,l\in\mathbb{Z}$ such that $$\cos[\tau(j+\dfrac{1}{2})\pi]>0, \cos[\tau(l+\dfrac{1}{2})\pi]<0.$$
\end{lemma}
Proof
We only prove there exists a suitable integer number $j\in\mathbb{Z}$ such that $\cos[\tau(j+\dfrac{1}{2})\pi]>0$. The existence of $l\in\mathbb{Z}$ is similar.
Take a nonnegative continuous periodic function $f$ with periodic $2$ with its support on interval $(2m-\dfrac{\tau}{2},2m+\dfrac{1-\tau}{2})$ for integer $m$.

By ergodic theorem \cite{probability}, we obtain that
$
\lim_{N\to+\infty}\dfrac{1}{N}\sum_{k=1}^N f(k\tau)=\dfrac{1}{2}\int_0^{2} f(x){\rm d}x>0.
$
This further indicates that there exist $j, m\in\mathbb{N}^*$ such that
$
2m-\dfrac{\tau}{2}<j\tau<2m+\dfrac{1-\tau}{2},
$
which implies that
$
2m\pi<\tau(j+\dfrac{1}{2})\pi<(2m+\dfrac{1}{2})\pi,
$
and therefore leads to that $\cos[\tau(j+\dfrac{1}{2})\pi]>0.$

\begin{lemma}\label{lemma13}
For $\lambda={\rm i}\beta$, $|\beta|\leq s_\epsilon\pi$, $\delta\in[-|\epsilon|,|\epsilon|]$, $|\epsilon|<\dfrac{1-\tilde{c}}{2}$, we have $h_\delta(\lambda)\neq c.$
\end{lemma}
Proof
By taking $\lambda={\rm i}\beta$, we compute that
\begin{small}$$
h_\delta(\lambda)=-\dfrac{1}{2}\left\{\cos [(2n+\delta)\beta]+\cos[(2n-2+\delta)\beta]\right\}
-\dfrac{\rm i}{2}\left\{\sin [(2n+\delta)\beta]+\sin[(2n-2+\delta)\beta]\right\}.
$$
\end{small}
If $h_\delta(\lambda)\in\mathbb{R}$, then we obtain that $\sin (2n+\delta)\beta+\sin(2n-2+\delta)\beta=0$, which, by Trigonometric Identities Equations, implies that $\cos\beta=0$ or $\sin[(2n-1+\delta)\beta]=0$.

If $\cos\beta=0$, we compute that $h_\delta(\lambda)=0\neq c$.
Now we turn to investigate the situation  $\sin[(2n-1+\delta)\beta]=0$. In this case, we get  that $(2n-1+\delta)\beta=k\pi, k\in\mathbb{Z}$. Furthermore, we compute that $h_\delta(\lambda)=(-1)^{k+1}\cos(\dfrac{k\pi}{2n-1+\delta})$. Now we prove that $h_\delta(\lambda)\neq c$. If not, by $|h_\delta(\lambda)|=|c|$, we get
$
|\cos(\dfrac{k\pi}{2n-1+\delta})|=\sin [\dfrac{\tilde{c}\pi}{2(2n-1)}].
$
By  Trigonometric Identities Equations, we obtain that
$\dfrac{k\pi}{2n-1+\delta}+ \dfrac{\tilde{c}\pi}{2(2n-1)}=(l+\dfrac{1}{2})\pi$,or~$\dfrac{k\pi}{2n-1+\delta}-\dfrac{\tilde{c}\pi}{2(2n-1)}=(l+\dfrac{1}{2})\pi$,
for some $l\in\mathbb{Z}$. Multiplying $\dfrac{2(2n-1)}{\pi}$ of both sides leads to that
$-\dfrac{2\delta k}{2n-1+\delta}+\tilde{c}=(2l+1)(2n-1)-2k$,~or $-\dfrac{2\delta k}{2n-1+\delta}-\tilde{c}=(2l+1)(2n-1)-2k.$
By using $(2n-1+\delta)\beta=k\pi$, we obtain $-\dfrac{2\delta\beta}{\pi}\pm\tilde{c}=(2l+1)(2n-1)-2k.$
However, by using $|\delta|\leq |\epsilon|$ and $|\beta|\leq s_\epsilon\pi$, we estimate that
$$|-\dfrac{2\delta\beta}{\pi}\pm\tilde{c}|\leq |\dfrac{2\delta\beta}{\pi}|+|\tilde{c}|\leq 2|\epsilon| s_\epsilon+\tilde{c}\leq |\epsilon|+\dfrac{C_1}{\pi}+\tilde{c}\leq  |\epsilon|+\dfrac{1-\tilde{c}}{2}+\tilde{c}<1.
$$
This is a contradiction because $(2l+1)(2n-1)-2k$ must be an odd number. This completes the proof.

\begin{lemma}\label{lemma14}
$\mathscr{F}_{\epsilon}$ is nonempty.
\end{lemma}
ProofSimilar to the proof of Lemma \ref{lemma13}, by taking $\lambda={\rm i}\beta$  and $h_\delta(\lambda)$ is a nonzero real number, we get that
$(2n-1+\delta)\beta=k\pi, k\in\mathbb{Z}$, $h_\delta(\lambda)=(-1)^{k+1}\cos(\dfrac{k\pi}{2n-1+\delta})$.

Condition \ref{condition} implies that $c=(-1)^n\sin[\dfrac{\tilde{c}\pi}{2(2n-1)}]$. By taking $h_\delta(\lambda)=c$, we obtain that
$\cos(\dfrac{k\pi}{2n-1+\delta})=(-1)^{n+k+1}\sin[\dfrac{\tilde{c}\pi}{2(2n-1)}].$
By  Trigonometric Identities Equations, we obtain that
$\dfrac{\pi}{2}-\dfrac{k\pi}{2n-1+\delta}=(n+k+1)\pi+\dfrac{\tilde{c}\pi}{2(2n-1)}+2l\pi$,
or
$
\dfrac{\pi}{2}-\dfrac{k\pi}{2n-1+\delta}+(n+k+1)\pi+\dfrac{\tilde{c}\pi}{2(2n-1)}=2l\pi+\pi,
$
for some $l\in\mathbb{Z}$.
Multiplying $\dfrac{2(2n-1)}{\pi}$ of both sides leads to that
$\dfrac{2\delta k}{2n-1+\delta}-\tilde{c}=4n^2+4kn+8ln-4l-1,$
$\dfrac{2\delta k}{2n-1+\delta}+\tilde{c}=-4n^2-4kn+8ln+4k-4l+1$.
By using $(2n-1+\delta)\beta=k\pi$, we obtain that
$\dfrac{2\delta\beta}{\pi}-\tilde{c}=4n^2+4kn+8ln-4l-1$, or
$\dfrac{2\delta\beta}{\pi}+\tilde{c}=-4n^2-4kn+8ln+4k-4l+1$.
By using $|\delta|\leq |\epsilon|$ and $|\beta|\leq S_\epsilon\pi$, we estimate that
$|\dfrac{2\delta\beta}{\pi}\pm\tilde{c}|\leq |\dfrac{2\delta\beta}{\pi}|+\tilde{c}\leq 2\epsilon S_\epsilon+\tilde{c}\leq 1+\tilde{c}$.
Note that $4n^2+4kn+8ln-4l-1$  is an integer number with the formation $4s-1$  for some  $s\in\mathbb{Z}$, thus we obtain that $\dfrac{2\delta\beta}{\pi}-\tilde{c}=-1=4n^2+4kn+8ln-4l-1$. Since $0<\tilde{c}<1$, we obtain that $\delta<0$ and $k<0$.
Similarly for the case $\dfrac{2\delta\beta}{\pi}+\tilde{c}=-4n^2-4kn+8ln+4k-4l+1$,
we obtain $\dfrac{2\delta\beta}{\pi}+\tilde{c}=1=-4n^2-4kn+8ln+4k-4l+1, \delta>0, k>0.$
Now we prove this Lemma by finding suitable $k,l\in\mathbb{Z}, \delta>0,\beta$. We consider the situation $k>0$. $-1=4n^2+4kn+8ln-4l-1$ is equivalent to
$k(1-n)+l(2n-1)=-n^2.$
Denote by $q$ as the largest integer number smaller than $(2n-1)S_\epsilon$. Since $1-n$ and $2n-1$ are coprime, there exist two integers $k^*,l^*$ with $k^*\in\{q-2n+2,q-2n+1,\cdots, q\}$ such that
$k^*(1-n)+l^*(2n-1)=-n^2$. We know that $(2n-1)S_\epsilon-k^*<2n-1$. Then by considering $\dfrac{2\delta^* k^*}{2n-1+\delta^*}=2\tilde{c}$, we choose $\delta^*=\dfrac{(2n-1)(1-\tilde{c})}{2k^*-1+\tilde{c}}$, and
$\beta^*=\dfrac{k^*\pi}{2n-1+\delta^*}$. Finally, we need to verify that $\delta^*\leq \epsilon$ and $\beta^*\leq S_\epsilon\pi$. By using $k^*>(2n-1)(S_\epsilon-1)$, we obtain that
$
\delta^*<\dfrac{(2n-1)(1-\tilde{c})}{2(2n-1)(S_\epsilon-1)-1+\tilde{c}}<\epsilon,
$
$
\beta^*<\dfrac{k^*\pi}{2n-1}<S_\epsilon\pi.
$
This completes the proof.

\begin{lemma}\label{lemma15}
$H_\epsilon({\delta})$ is increasing with respect to $\delta\in(0,\epsilon]$.
\end{lemma}

ProofSimilar to the proof in Section \ref{section3}, we prove that for each element $\delta_*\in\mathscr{F}_{\epsilon}$, we have $h_{\delta^*}({\rm i}\beta_*)=c$. Then there exists an implicit function $\lambda(\delta)$ such that $h_{\delta}(\lambda(\delta))=c$, $\lambda(\delta_*)={\rm i}\beta_*$. Furthermore, we have
${\rm Sgn}[{\rm Re}\lambda'(\delta_*)]>0$.

By taking derivatives of $\delta$ of $h_{\delta^*}({\rm i}\beta_*)=c$, we obtain that
$\lambda'(\delta)=-\dfrac{\dfrac{\partial h_\delta}{\partial \delta}}{h'_{\delta}(\lambda)}=-\dfrac{\lambda(1+{\rm e}^{-2\lambda})}{(2n+\delta)(1+{\rm e}^{-2\lambda})-2{\rm e}^{-2\lambda}}.$
By taking $\lambda={\rm i}\beta_*$, we compute
$
{\rm Sgn}[{\rm Re}\lambda'(\delta_*)]=-{\rm Sgn}[\beta_*\sin(2\beta_*)].
$
The first case is $\beta_*>0$, then there exists integer number $k_*>0,n,l_*$ such that
$\dfrac{\pi}{2}-\dfrac{k_*\pi}{2n_*-1+\delta}+(n_*+k_*+1)\pi+\dfrac{\tilde{c}\pi}{2(2n_*-1)}=2l_*\pi+\pi$. Multiplying $2$ on both sides leads to that
$
2\beta_*=\dfrac{\tilde{c}\pi}{2n-1}+2(n+k_*+1)\pi-4l_*\pi-2\pi+\pi,
$
which implies that
$\sin(2\beta_*)=-\sin [\dfrac{\tilde{c}\pi}{2n-1} ]<0$. Thus we have  ${\rm Sgn}[{\rm Re}\lambda'(\delta_*)]>0$. For the case $\beta^*<0$, the proof is similar. This completes the proof.

\begin{lemma}\label{lemma16}
For $\epsilon>0$ small efficiently,  Eq. $\eqref{202302042155}$ has at least one root located in $\Big\{\lambda\in\mathbb{C}_0\cup\mathbb{C}_+\Big||{\rm Im}\lambda|<S_\epsilon\pi\Big\}$ when we take $c=c>0$, $\tau=\epsilon$. Here $S_\epsilon$ is the smallest integer number such that $S_\epsilon>\dfrac{1}{\epsilon}$.\end{lemma}
Proof
The idea of the proof is quite similar to Lemma \ref{lemma7}, \ref{lemma8} and \ref{lemma11}. We use the notation $M_{-S_{\epsilon},S_{\epsilon}}(c)$ to denote the number of root of  Eq. $\eqref{202302042155}$ located in $\Big\{\lambda\in\mathbb{C}_0\cup\mathbb{C}_+\Big||{\rm Im}\lambda|<S_\epsilon\pi\Big\}$. In this notation, we fix $\tau=\epsilon$ and allow $c$ to vary on the whole interval $(0,+\infty)$. If we denote by a set
$$
\mathscr{H}\triangleq\left\{g(\lambda) \Bigg| \lambda\in\mathbb{C}_0, |{\rm Im}\lambda|\leq S_\epsilon \pi \right\} \cap \mathbb{R},
$$
where $$g(\lambda)\triangleq-\dfrac{1}{2}[{\rm e}^{\epsilon\lambda}+{\rm e}^{(\epsilon-2)\lambda}].$$
Using the same method as in the poof of Lemma \ref{lemma7} and \ref{lemma11}, we obtain that $M_{-S_{\epsilon},S_{\epsilon}}(c)$ is a constant on $(1,+\infty)$ and $M_{-S_{\epsilon},S_{\epsilon}}(c)$ is monotonically decreasing on $\mathbb{R}_+\backslash\mathscr{H}$. Finally, we only to validate that $M_{-S_{\epsilon},S_{\epsilon}}(c)$ is nonzero on  $(1,+\infty)$. We use the same method as in the proof of Lemma \ref{lemma8}. We prove that there exists $\lambda\in\Big\{\lambda\in\mathbb{C}_0\cup\mathbb{C}_+\Big||{\rm Im}\lambda|<S_\epsilon\pi\Big\}$ such that $g(\lambda)\in(1,+\infty)$.

Suppose that $\lambda=p+{\rm i}q, p>0,|q|<S_\epsilon\pi$. Thus, we get that
$$
g(\lambda)=-\dfrac{1}{2}[{\rm e}^{\epsilon p}\cos(\epsilon q)+{\rm e}^{\epsilon p-2p}\cos(\epsilon q-2q)]-\dfrac{{\rm i}}{2}[{\rm e}^{\epsilon p}\sin(\epsilon q)+{\rm e}^{\epsilon p-2p}\sin(\epsilon q-2q)].
$$

Denote by $q^*\triangleq\dfrac{\pi}{\epsilon}<S_{\epsilon}\pi$, there are two different situations. Case I: $\sin(\epsilon q^*-2q^*)=0$. Case II: $\sin(\epsilon q^*-2q^*)\neq 0$.

For Case I, when $q=q^*$, we have ${\rm Im}~g(\lambda)=0$, $\cos(\epsilon q^*)=-1$, $\cos(\epsilon q^*-2q^*)=\pm 1$. Then $g(\lambda)=\dfrac{1}{2}{\rm e}^{\epsilon p}(1\pm {\rm e}^{-2p})$. For sufficiently large $p>0$, we have $g(\lambda)\in(1,+\infty)$.

For Case II, considering that ${\rm e}^{\epsilon p}\sin(\epsilon q)+{\rm e}^{\epsilon p-2p}\sin(\epsilon q-2q)=0$ yields that
${\rm e}^{2p}=-\dfrac{\sin(\epsilon q-2q)}{\sin(\epsilon q)},g(\lambda)=\dfrac{{\rm e}^{\epsilon p}\sin(2q)}{2\sin(\epsilon q-2q)}$. Since $\sin(\epsilon q^*)=0, \sin(\epsilon q^*-2q^*)\neq 0$,
we can take a sequence of $\{p_n\},\{q_n\}$ such that $\lim_{n\to+\infty}q_n=q^*,\lim_{n\to+\infty}p_n=+\infty$ while   ${\rm e}^{2p_n}=-\dfrac{\sin(\epsilon q_n-2q_n)}{\sin(\epsilon q_n)}.$ Thus, we have $\lim_{n\to+\infty}\dfrac{\sin(2q_n)}{\sin(\epsilon q_n-2q_n)}=1$. If we take $\lambda_n=p_n+{\rm i}q_n$, we obtain that $\lim_{n\to+\infty}g(\lambda_n)=+\infty$ and thus complete the proof.

\end{document}